\documentclass[a4paper,12pts]{amsart}
\usepackage{graphicx} 
\usepackage{amsthm}
\usepackage{amssymb}
\usepackage{amsfonts}
\usepackage{mathtools}
\usepackage{mathrsfs}
\usepackage{graphicx}
\usepackage{color}
\usepackage{enumitem}
\usepackage{setspace}
\usepackage[english]{babel}
\usepackage{csquotes}
\usepackage{etoolbox}
\usepackage{thmtools}
\usepackage{yfonts}
\usepackage{mfirstuc}
\usepackage{stix}
\usepackage{fancyvrb}
\usepackage{textgreek} 

\usepackage{amsmath}
\usepackage{amsthm, amssymb, amsfonts}
\usepackage{graphicx, color}
\usepackage{caption, subcaption} 
\usepackage{multicol} 
\usepackage{stmaryrd}
\usepackage{tikz-cd} 
\usepackage{lipsum}  
\usepackage[capitalise,noabbrev, nameinlink]{cleveref} 
\Crefname{part}{Part}{Parts}
\Crefname{step}{Step}{Steps}
\Crefname{prop}{Proposition}{Propositions}
\Crefname{prob}{Problem}{Problems}
\newcommand{\forces}[1]{\Vdash\text{\say{#1}}}
\newcommand{\Lim}{\text{Lim}}
\newcommand{\Fin}{\text{Fin}}

\usepackage{enumitem}
\renewcommand{\leq}{\leqslant}
\renewcommand{\geq}{\geqslant}

\renewcommand{\epsilon}{\varepsilon}

\newcommand{\stickT}{%
\setbox255=\hbox{\raise1ex\hbox{$\hspace{0.2pt}\,\bullet\,$}}
\mathord{\rlap{\hbox to\wd255{\hss\hbox{$|$}\hss}}
\box255}
}
\newcommand{\stickS}{%
\setbox255=\hbox{\raise0.6ex\hbox{$\scriptstyle\bullet$}}
\mathord{\rlap{\hbox to\wd255{\hss\hbox{$\scriptstyle|$}\hss}}
\box255}
}
\newcommand{\stick}{{\mathchoice{\stickT}{\stickT}{\stickS}{\stickS}}}

\numberwithin{equation}{section}
\theoremstyle{plain}
\newtheorem{theorem}[equation]{Theorem}
\newtheorem{lemma}[equation]{Lemma}
\newtheorem{proposition}[equation]{Proposition}
\newtheorem{corollary}[equation]{Corollary}

\theoremstyle{definition}
\newtheorem{definition}[equation]{Definition}

\newtheorem{problem}[equation]{Problem}

\theoremstyle{remark}

\newtheorem{rem}[equation]{Remark}

\usepackage{mathtools}
\usepackage{amsfonts}
\usepackage{amssymb}
\usepackage{amsthm}
\usepackage{amsmath}
\usepackage{dirtytalk}
\usepackage{mathrsfs}  
\usepackage{lineno}
\newtheorem{remark}{Remark}
 \usepackage{multicol}  
 \usepackage{bbm}
 \usepackage{tabularx}
\newenvironment{claimproof}[1][\proofname]
{\begin{proof}[#1]}
{\end{proof}}

\title{Multiple gaps and some finitizations of Club and CH}
\author{Jorge Cruz}

\begin{document}
\maketitle
\begin{abstract}
We continue the development of the theory of capturing schemes over $\omega_1$ by analyzing the relation between the capturing construction schemes (whose existence is implied by Jensen's $\Diamond$-principle) and both the Continuum Hypothesis and Ostaszewski's $\varclubsuit$-principle. Formally, we show that the property of being capturing can be viewed as the conjunction of two properties, one of which is implied by $\varclubsuit$ and the other one by CH.  We  apply these principles to construct multiple gaps, entangled sets and metric spaces without uncountable monotone subspaces.
\end{abstract}
A a tuple $\mathfrak{L}=(\mathcal{L}_i\,:\,i<n)$ of families of subsets of a set $X$ is said to form a \emph{pregap} if $\mathcal{L}_i$ is orthogonal to $\mathcal{L}_j$ for any two distinct $i$ and $j$. That is, $L\cap R=^*\emptyset$ for all $L\in \mathcal{L}_i$ and $R\in \mathcal{L}_j$.  A function  $c:n\longrightarrow \mathscr{P}(X)$ is said to \emph{separate} $\mathfrak{L}$ in case $L\subseteq^*c(i)$ and $C\cap R=^*\emptyset$ for each $i<n$ and $L\in\mathcal{L}_i$. The pregap $\mathfrak{L}$ is said to be an $n$-\textit{gap} (or simply, multiple gap) if there is no function $c$ separating it. 
Given a family $\langle (X_i,\leq_i)\rangle_{i<n}$  of partial orders, we say that  $\mathfrak{L}$ is an \emph{$(X_i\,:\,i<n)$-pregap (resp. $(X_i\,:\,i<n)$-gap)} if $(\mathcal{L}_i,\subseteq^*)$ is isomorphic to $(X_i,\leq_i)$ for each $i<n$.  If moreover $(X_i,\leq_i)=(X,<)$ for any $i$, we simply refer to it as an \emph{$X$-$n$-pregap} or $X$-$n$-gap respectively.\\
In 1909, Hausdorff constructed for the first time an $(\omega_1,\omega_1)$-gap without appealing to any extra axioms. This clever idea became a blueprint for the construction of several uncountable objects in the following years. 
Motivated by Amir's result (\cite{Amir}) stating that the space $l_\infty/c_0$ is not injective, Avil\'es and Todor\v{c}evi\'c introduced in \cite{multiplegaps} the concept of an $n$-gap as a natural generalization of the well-studied notion of ($2$-)gaps presented by Hausdorff. In there, they showed that unlike usual gaps, Hausdorff's result about the existence of $(\omega_1,\omega_1)$-gaps in $ZFC$ can not be extended to greater dimensions without assuming extra axioms. Precisely, they showed that under $MA_\theta(\sigma-k\text{-linked})$ there are no $n$-gaps formed by $\theta$-generated ideals in $\mathscr{P}(X)$ for $n>k$ and $X$ countable. Particularly, this means that the non-existence of $\omega_1$-$n$-gaps is consistent with the usual axioms of set theory for each $n\geq 3$.  On the other hand, they proved that CH implies the existence  of $\omega_1$-$n$-gaps for every $n$. The theory  has been further developed in papers such as \cite{analyticstrongngaps}, \cite{finitebasismultiplegaps} and \cite{isolatingsubgaps}, with a particular emphasis on the study of  multiple gaps formed by analytic families. Beyond the results mentioned above, not as much is known about multiple gaps generated by families of size $\omega_1$. \\ In this manuscript we will develop powerful combinatorial tools  which will allow us to further expand our knowledge of the structure of $n$-gaps over $\omega$. Not only that, but we will use these same tools to:
\begin{itemize}
\item Show the existence of entangled sets of reals under the new axiom $FCA^\Delta$, weaker than CH, compatible with arbitrarily large continuum and which holds, for example, in the Cohen model. 
\item Present some partial progress on a question of Hru\v{s}\'ak and Zindulka in \cite{hrusakzindulka}. Specifically, we will prove that under $CA^\Delta$, a weaker version of $FCA^\Delta$, there exists an uncountable metric space without uncountable monotone subspaces.

\end{itemize}
In \emph{Section \ref{sectiongaps}}, we define the concept of an $n$-inseparable AD family and the AD representation of a partial order. We state Theorem \ref{omegainseptheorem} and use it to show that both under CH and $\varclubsuit$ there are (among other types of multiple gaps) $\omega_1$-$n$-gaps over $\omega$. In particular, this means that $\omega_1$-$n$-gaps may exist in the Sacks model. We then generalize the concept of the gap cohomology group introduced by Talayco in \cite{cohomologytalayco} and  use AD representations to show the existence of many structures with big gap cohomology groups  (see Theorem \ref{biggapcohomologytheorem}).  In \emph{Section \ref{sectionfinitizationcapturingaxioms}}, we recall the theory of construction schemes introduced by Todor\v{c}evi\'c in \cite{schemenonseparablestructures} and define the capturing axioms $CA^\rho$, $CA^\Delta$ (a variant of $FCA^\Delta$) and some of their variations.  In \emph{Section \ref{sectionproofomegainseptheorem}}, we prove Theorem \ref{omegainseptheorem} by showing that it follows from both $CA^\rho$ and $CA^\Delta$. In \emph{Section \ref{sectionentangled}}, we provide two more applications of the axiom $FCA^\Delta$ by constructing entangled sets of reals and an uncountable metric space without uncountable monotone subspaces. In \emph{Section \ref{sectionforcing}}, we provide a quick review of the forcing $\mathbb{P}(\mathcal{F})$. This forcing is the main tool needed to build construction schemes in a recursive manner. In \emph{Section \ref{sectionclub}}, we show that $CA^\rho$ follows from the  $\varclubsuit$-principle and in \emph{Section \ref{sectionch}} we show that $FCA^\Delta$ follows from CH. Finally, in \emph{Section \ref{sectionfinal}} we pose some open problems.\\\\
The author would like to thank Osvaldo Guzm\'an, Michael Hru\v{s}\'ak and Stevo Todor\v{c}evi\'c for their valuable comments and support during the preparation of this work. Part of it is based on a previous version of \cite{schemescruz} which is a joint work with Guzm\'an and Todor\v{c}evi\'c. 
\section{Notation}
The notation and terminology used here mostly standard and it follows \cite{schemescruz}.  Given a set $X$ and a (possibly finite) cardinal $\kappa$, $[X]^\kappa$ denotes the family of all subsets of $X$ of cardinality $\kappa$. The sets $[X]^{<\kappa}$ and $[X]^{\leq \kappa}$ have the expected meanings. The family of all non-empty finite sets of $X$ is denoted by $\Fin(X)$. That is, $\Fin(X)=[X]^{<\omega}\backslash\{\emptyset\}$. $\mathscr{P}(X)$ denotes the power set of $X$. By $\text{Lim}$  we mean the  set of limit ordinals strictly smaller than $\omega_1$. Given sets $X$ and $Y$ two sets of ordinals, we write $X<Y$ whenever $\max(X)<\min(Y)$ or $X=\emptyset$. For a set of ordinals $X$, we denote by $ot(X)$ its order type. We identify $X$ with the unique strictly increasing function $h:ot(X)\longrightarrow X$. In this way, $X(\alpha)=h(\alpha)$ denotes the $\alpha$ element of $X$ with respect to its increasing enumeration. Analogously, $X[A]=\{X(\alpha)\,:\,\alpha\in A\}$ for $A\subseteq ot(X)$.  A family $\mathcal{D}$ is called a \textit{$\Delta$-system} with root $R$ if $|\mathcal{D}|\geq 2$ and $X\cap Y= R$ whenever $X,Y\in \mathcal{D}$ are different.  If moreover, $R<X\backslash R$ for any and $X\backslash R<Y\backslash R$ or viseversa for any two $X,Y\in \mathcal{D}$, we call $\mathcal{D}$ a root-tail-tail $\Delta$-sytem.
\section{Almost disjoint representations and multiple gaps}\label{sectiongaps}
Recall that an \emph{almost disjoint family} (shortly, AD family) over a set $X$,\footnote{If $X $ is not mentioned, we will assume that it is $\omega$.} is a family $\mathcal{A}$ of infinite subsets of $X$ so that $A\cap B=^*\emptyset$ for any two distinct $A,B\in \mathcal{A}$.

 We say that $\mathcal{A}$ is:
\begin{itemize}
    \item $n$-\emph{inseparable} if for any  pairwise disjoint $\mathcal{L}_0,\dots,\mathcal{L}_{n-1}\in [\mathcal{A}]^{\omega_1}$, the tuple $(\mathcal{L}_i\,:\,i<n)$ forms an $n$-gap. We suppress the index if $n=2$.
    \item \emph{$\omega$-inseparable} if $\mathcal{A}$ is $n$-inseparable for each $n\geq 2$.
    \item  \emph{Jones} if for any two disjoint  $\mathcal{L}\in [\mathcal{A}]^{\omega}$ and $\mathcal{R}\in [\mathcal{A}]^{\omega_1}$,  $(\mathcal{L},\mathcal{R})$ is not a gap. 
  \end{itemize}
 The first construction of an inseparable family was done in \cite{luzin1947subsets} by Luzin. A Jones family was implicitly constructed by F. B. Jones in \cite{jones1937concerning}. At first glance,  inseparable and Jones seem to be incompatible properties. However,  a construction of a inseparable-Jones family was obtained in \cite{guzman2019mathbb} by Guzm\'an, Hru\v{s}\'ak and Koszmider (building from work by Koszmider in \cite{onconstructionswith2cardinals}). \\
 The classical constructions of both an inseparable AD family  and an $(\omega_1,\omega_1)$-gap share the same basic idea. However, the reader might agree that from these two, the latter is much  harder to build. Having this in mind, one might ask whether it is possible to build an $(\omega_1,\omega_1)$-gap in an easier way by using the elements of an already existing inseparable family as building blocks. A natural attempt of doing this is by starting with an almost disjoint family $\mathcal{A}=\{A_\alpha\,:\,\alpha\in \omega_1\}$ and build a pregap $(\mathcal{L},\mathcal{R})=(L_\alpha,R_\alpha)_{\alpha\in \omega_1}$ as follows: We start by defining $L_0=A_0$ and $R_0=A_1$. If $\alpha=\gamma+n$ for some limit natural $n$ greater than $0$, we let $L_\alpha:=L_{\alpha-1}\cup A_{\gamma+2n}$ and $R_\alpha:=R_{\alpha-1}\cup A_{\gamma+2n+1}$. Assume that at limit steps $\gamma$ we are able to find  \begin{center}$(\star)$
 \\two disjoint sets $L$ and $R$ so that $L_\alpha\subseteq^*L$, $R_\alpha\subseteq^* R$,  and $A_{\beta+2n+1}\cap L=^*\emptyset=^*A_{\beta+2n}\cap R$ for every $\alpha<\gamma$ and each limit ordinal $\beta$.\end{center} Then we can continue the construction by defining $L_\gamma:=L\cup A_\gamma$ and $R_\gamma:=R\cup A_{\gamma+1}$. At the end, we will have built a pregap $(\mathcal{L},\mathcal{R})=(L_\alpha,R_\alpha)_{\alpha\in \omega_1}$. If $\mathcal{A}$ is inseparable, it is not hard to see that $(\mathcal{L},\mathcal{R})$ will in fact be a gap. 
So the question is: Is there an inseparable almost disjoint $\mathcal{A}$ for family for which $(\star)$ holds at every limit step? One can prove that the answer to this question is \say{yes} under $\mathfrak{b}>\omega_1$. In fact, under this hypothesis it is possible to carry the construction for any inseparable AD family. Another trivial observation is that any such $\mathcal{A}$ should at least be Jones. For this same reason, we can not expect to show that the construction described above will work for every $\mathcal{A}$. This is because while it is true that there exists an inseparable Jones family in $ZFC$, under $\mathfrak{b}=\omega_1$ there is also one inseparable family which is not Jones.\\
The next definition is based on the discussion  presented above. Informally, it captures the idea of being able to embed a given partial order in the power set of a set $N$ using an almost disjoint family over $N$ as
\say{building blocks} of this embedding.

\begin{definition}[Almost disjoint representation]\label{Luzinrepdef} Let $(X,\leq)$ be a partial order of size $\omega_1$. An \textit{almost disjoint representation of $X$} (Shortly, AD-representation) is an ordered pair $(\mathcal{T},\mathcal{A})$ consisting of two families of infinite subsets of $\omega$ indexed as $\langle T_x\rangle_{x\in X}$ and $\langle A_x\rangle_{x\in X}$ respectively. Furthermore, $\mathcal{A}$ is an AD family and for all $x,y\in X$ the following properties hold:
\begin{enumerate}[label=$(\alph*)$]
\item $A_x\subseteq T_x.$

\item If $y\not\leq x$, then $A_y\subseteq^*T_y\backslash T_x$.
\item If $\inf(x,y)$ exists, then $T_x\cap T_y=^*T_{\inf(x,y)}.$

\item If $x$ is succesor-like, then $$T_x\backslash\big(\bigcup\limits_{z\in pred(x)}T_z\big)=^*A_x.$$
\item If $x$ and $y$ are incompatible, then $T_x\cap T_y=^*\emptyset.$

\end{enumerate}
In particular, $(\mathcal{T},\subseteq^*)$ is order isomorphic to $X$ by virtue of the points $(b)$ and $(c)$.  In the case there is $(\mathcal{T},\mathcal{A})$ an AD representation of $X$ we  will say that $\mathcal{A}$ \textit{codes} $X$.
\end{definition}
The construction of a gap $(\mathcal{L}, \mathcal{R})$ via an AD family $\mathcal{A}$ described at the beginning of this section would actually yield an AD representation $(\mathcal{L}\cup \mathcal{R},\mathcal{A})$ of the partial order $(2\times \omega_1,\leq)$ where $(i,\alpha)\leq (j,\beta)$ if and only if $i=j$ and $\alpha\leq \beta$. So now the question is: Is there an inseparable AD family $\mathcal{A}$ which which codes $2\times \omega_1$? We will answer this question in the affirmative not only for $2\times \omega_1$ but, up to cofinal suborders, for any possible partial order of size $\omega_1$.\\
Let us say that a partial order $(X,\leq)$ is \emph{$\omega_1$-like} if it is well-founded, $|X|=\omega_1$ and $|(-\infty,x)|\leq \omega$ for all $x\in X$.

\begin{proposition}\label{propositionwellfounded}Let $(Y,\leq)$ be a partial order of cardinality $\omega_1$. There is a well-founded cofinal $X\subseteq Y$ with $|(-\infty, x)_X|\leq \omega$ for each $x\in X$.
\begin{proof} Enumerate $Y$ as $\langle y_\alpha\rangle_{\alpha\in \omega_1}$.  Define $X$ as the set of all $y_\beta$'s such that $y_\beta\not\leq y_\alpha$ for each $\alpha<\beta$. Of course $X$ is well-founded as for each $x,z\in X$ with $x=y_\alpha$ and $z=y_\beta$ it happens that if $x<z$ then $\alpha<\beta$. Because of this we also have that $(-\infty,\, y_\beta)_X\subseteq \{y_\alpha\,:\,\alpha<\beta\}$ for each $y_\beta\in X$. Lastly $X$ is cofinal in $Y$ because for each $\beta\in \omega_1$ the element $y_{\xi_\beta}\in X$ where  $\xi_\beta=\min(\,\alpha\leq \beta\,:\, y_\alpha \geq y_\beta\,)$.
\end{proof}
\end{proposition}
We are ready to present the two main theorems regarding almost disjoint families and multiple gaps.
\begin{theorem} There is an inseparable AD family $\mathcal{A}$ which codes any $\omega_1$-like order.
\end{theorem}
\begin{theorem}[Under CH or $\varclubsuit$-principle]\label{omegainseptheorem} There is an $\omega$-inseparable AD family $\mathcal{A}$ which codes any $\omega_1$-like order.
\end{theorem}
Theorem  \ref{omegainseptheorem} can be stated in a more general way. We will present two new axioms, namely $CA^\Delta$ and $CA^\rho$, which are consequences of CH and $\varclubsuit$ respectively. These two principles are compatible with arbitrarily large values of the continuum. In particular, they both hold in the Cohen model. Actually, we will show that the existence of an AD family satisfying the conclusion Theorem \ref{omegainseptheorem} follows from both $CA^\Delta$ and $CA^\rho$.

The existence of an $\omega$-inseparable AD family is already interesting as it is the natural generalization to the realm of $n$-gaps of the well-studied inseparable families. As we show in the following corollary, Theorem \ref{omegainseptheorem} yields new situations in which $\omega_1$-n-gaps exist.
\begin{corollary}Assume there is an $n$-inseparable AD family $\mathcal{A}$ which codes any $\omega_1$-like order and let $\langle (X_i,<_i)\rangle_{i<n}$ is a family of partial orders of cofinality $\omega_1$. Then there are cofinal $X'_i\subseteq X_i$ for each $i$, so that there is an $(X'_i\,:\,i<n)$-gap.
\begin{proof}
 Without loss of generality we may assume that $X_i\cap X_j=\emptyset$ for any two distinct $i$ and $j$. By applying Proposition \ref{propositionwellfounded} we get, for each $i$, an $\omega_1$-like  $X'_i\subseteq X_i$ an cofinal subsets of $X_i$. Let $Z=\bigcup\limits_{i<n}X'_i$ 
and $<_Z:=\big(\bigcup\limits_{i<n} <_i\big)\cap Z^2$ . That is, $x<_Z y$ if and only if there is $i$ for which $x,y\in X'_i$ and $x<_{X_i} y$. It is straightforward that $(Z,<_Z)$ is an $\omega_1$-like order. 
Furthermore, if $x\in X'_i$ and $y\in X_j'$  for $i\not=j$, then there is no $z\in Z$ with $z\leq_Z x, y.$ Let $(\mathcal{T},\mathcal{A})$ be an almost disjoint representation of $Z$ wifh $\mathcal{A}$ being $n$-inseparable. Due to  Definition \ref{Luzinrepdef} it follows that $(\langle T_x\rangle_{x\in X'_i}\,:\,i<n)$  is an $(X'_i\,:\,i<n)$-pregap.
We claim that this is in fact a $n$-gap. For this just note that any function $c$ separating $(\langle T_x\rangle_{x\in X'_i}\,:\,i<n)$ would also separate $(\langle A_x\rangle_{x\in X'_i}\,:\,i<n)$ by the point (a) of Definition \ref{Luzinrepdef}. But this is impossible since $\mathcal{A}$ is $n$-inseparable and each $X'_i$ is uncountable.
\end{proof}
\end{corollary}
Particularly, we get the following corollary. As a direct consequence of it, we also have that there are $\omega_1$-n-gaps in the Sacks model. 
\begin{corollary}[Under $\varclubsuit$] There are $\omega_1$-n-gaps for every $2\leq n\in \omega$.
\end{corollary}
In order to keep illustrating the power of the concept of almost disjoint representations, we now extend the work done by Talayco in \cite{cohomologytalayco}. Specifically, we  study the size of groups defined in a similar context as Talayco and prove that, under certain assumptions, these groups are as big as possible.\\
We say that $(X,\leq,\wedge)$ is a \emph{lower semi-lattice} if $(X,\leq)$ is a partial order and for any $x,y\in X$, $\inf(x,y)$ exists and its equal to $x\wedge y$. A family $\mathcal{T}\subseteq \mathscr{P}(\omega)$ is said to be a \emph{$*$-lower semi-lattice} if its projection to the Boolean algebra $\mathscr{P}(\omega)/\Fin$ is a lower semi-lattice inside $(\mathscr{P}(\omega)/\Fin,\,\subseteq,\, \cap )$ That is, for any two $A,B\in \mathcal{T}$ there is $C\in \mathcal{T}$ so that $A\cap B=^*C$.
Note that towers are particular cases of $*$-lower semi-lattices.
\begin{definition}[Coherent Subsystems]\label{cohsystemsdef} Let $\mathcal{T}$ be a $*$-lower semi-lattice over a countable set $X$. A function $g:\mathcal{T}\longrightarrow \mathscr{P}(X)$ is said to be a \textit{coherent subsystem of $\mathcal{T}$} if for any $A,B\in \mathcal{T}$, the following happens:
\begin{enumerate}[label=$(\arabic*)$]
\item $g(A)\subseteq A$,
\item if $B\subseteq^* A$ then $g(A)\cap B=^*g(B)$.
\end{enumerate}
 We say that $g$ is \textit{trivial} if there is $C\in \mathscr{P}(X)$ so that $C\cap A=^*g(A)$ for any $A\in \mathcal{T}$. In this case, we say that $C$ \textit{trivializes} $g$. The set of all coherent subsystems of $\mathcal{T}$ is denoted as $C(\mathcal{T})$ and the set of all trivial coherent subsystem is denoted as $Tr(\mathcal{T}).$
\end{definition}
\begin{definition}Let $\mathcal{T}$ be a $*$-lower semilattice over a countable set $X$ and $g$ be a coherent subsystem of $\mathcal{T}$. We denote the sequence $(g(A),A\backslash g(A))_{A\in\mathcal{T}}$ as $\mathcal{C}(\mathcal{T},g).$
    
\end{definition}
The following lemma is easy. We prove it just to emphasise the reason for defining coherent subsystems just for $*$-lower semi-lattices.
\begin{lemma}\label{lemmacoherentsubsystempregap}Let $\mathcal{T}$ be a $*$-lower semi-lattice and let $g$ be a coherent subsystem of $\mathcal{T}$. Then $\mathcal{C}(\mathcal{T},g)$ is a pregap. Furthermore, $g$ is non-trivial if and only if $\mathcal{C}(\mathcal{T},g)$ forms a gap.

\begin{proof}Let $A,B\in \mathcal{T}$. We need to prove that $g(A)\cap (B\backslash g(B))=^*\emptyset.$ For this let $C\in \mathcal{T}$ be so that $C=^*A\cap B$. Then \begin{align*}g(A)\cap B=^* (g(A)\cap A)\cap B=^*g(A)\cap C
=^*g(C)=^*g(B)\cap C\subseteq ^*g(B).
\end{align*}
This finishes the argument. Now we will prove the proposed equivalence.
\begin{claimproof}[Proof of $\Rightarrow$] We will prove this by contra-positive. Suppose that there is a set $C$ separating the pregap $\mathcal{C}(\mathcal{T},g)$. For any $A\in \mathcal{T}$ we have that $g(A)\subseteq^* C$ and $C\cap (A\backslash g(A))=^*\emptyset.$ In this way, $C\cap A=^*g(A)$. We conclude that $C$ trivializes $g$. 
\end{claimproof}
\begin{claimproof}[Proof of $\Leftarrow$] Again, by contra-positive. Suppose that there is $C$ which trivializes $g$. Then $C\cap A=^*g(A)$ for any $A\in \mathcal{T}$. In particular $g(A)\subseteq^* C$ and $C\cap (A\backslash g(A))=\emptyset$. We conclude that $C$ separates $\mathcal{C}(\mathcal{T},g)$.
\
\end{claimproof}
\end{proof}
\end{lemma}

\begin{definition}Let $\mathcal{T}$ be a $*$-lower semi-lattice over a set $X$. For $f,g\in C(\mathcal{T})$ we can define $f\cdot g\in C(\mathcal{T})$ given as: $$(f\cdot g)(A)= f(A)\Delta g(A).$$
\end{definition}
It is easy to see that $C(\mathcal{T})$ is an abelian (Boolean)\footnote{A group $(G,\cdot, e)$ is Boolean if $g\cdot g=e$ for any $g\in G$. Equivalently, $G$ is a vector space over $\mathbb{Z}_2$. In particular, the isomorphism type of $G$ is completely determined by its cardinality.} group with respect to this operation and $Tr(\mathcal{T})$ is a (normal) subgroup of it.  The quotient of these two groups is what we call \emph{the gap cohomology group of $\mathcal{T}$}.
\begin{definition}[Gap cohomology group] We define the \emph{gap cohomology group of $\mathcal{T}$} as the quotient of these two groups. That is,  $$G(\mathcal{T})=C(\mathcal{T})/Tr(\mathcal{T}).$$
We  say that $f,g\in C(\mathcal{T})$ are \emph{cohomologous} if the class of $f$ is equal to the class of $g$ inside $G(\mathcal{T})$, or equivalently, if $f\cdot g\in Tr(\mathcal{T})$. 
\end{definition}
\begin{rem}If $\mathcal{T}$ has size $\omega_1$, then $|G(\mathcal{T})|\leq 2^{\omega_1}.$
    
\end{rem}
In \cite{cohomologytalayco}, Talayco defined the gap cohomology group for the particular case of $\omega_1$-towers. In there he showed that the $2^\omega\leq |G(\mathcal{T})|$ for any $\omega_1$-tower $\mathcal{T}$. It is worth noting that by that time, Todor\v{c}evi\'c had already proved (in the 1980's) that an $(\omega_1, \omega_1)$-gap can be inserted into any tower  (see \cite{topicssettheory} or \cite{Walksonordinals}).\\
Talayco's result is based on the following theorem. \footnote{In a private communication Todor\v{c}evi\'c has pointed that his $(\omega_1, \omega_1)$-gap construction inside a given tower can easily be modified to give this result as well.}
\begin{theorem}[The $\aleph_0$ gap theorem] let $\mathcal{T}=\langle T_\alpha\rangle_{\alpha\in \omega_1}$ be an $\omega_1$-tower. There is  a sequence of functions $\langle f_\alpha\rangle_{\alpha\in\omega_1}$ with the following properties:
\begin{itemize}
\item $\forall \alpha\in \omega_1\,(\,f_\alpha:T_\alpha\longrightarrow \omega\,)$,
\item $\forall \alpha<\beta\in \omega_1\, (\,f_\alpha|_{T_\alpha\cap T_\beta}=^*f_\beta|_{T_\alpha\cap T_\beta}\,),$
\item $\forall m<n\in\omega\,(\, (f_\alpha^{-1}[\{m\}],f_{\alpha}^{-1}[\{n\}]\,)_{\alpha\in\omega_1}\textit{ is an $(\omega_1,\omega_1)$-gap }).$ 

\end{itemize}
\end{theorem}
He also proved that the $\Diamond$-principle implies that the size of this group is always $2^{\omega_1}$. Later, Todor\v{c}evi\'{c} noted that the $\aleph_0$ gap theorem already gives you the previous conclusion without assuming any extra axioms (see Theorem 32 in \cite{cohomologytalaycotrees}). Namely, if $2^{\omega}=2^{\omega_1}$ the result is clear due to Talayco's calculations. On the other hand, if $2^{\omega}<2^{\omega_1}$ then the result follows from an easy counting argument involving the quotient structure of $G(\mathcal{T}).$
In \cite{Agapcohomologygroup}, Morgan constructed an $\omega_1$-tower whose gap cohomology group can be explicitly calculated. This construction was carried through the use of morasses. Finally, in \cite{coherentfamilyoffunctions}, Farah improved the $\aleph_0$ gap theorem. This was achieved by considering the following concept.
\begin{definition}[Coherent families of functions]\label{coherentfamiliesdef} A \textit{coherent family of functions} supported by an $\omega_1$-tower $\langle T_\alpha\rangle_{\alpha\in\omega_1\backslash \omega}$ is a family of functions $\langle f_\alpha \rangle_{\alpha\in\omega_1}$ such that:
\begin{enumerate}[label=$(\arabic*)$]
    \item $\forall \alpha\in \omega_1\backslash \omega\,(\,f_\alpha:T_\alpha\longrightarrow \alpha\,),$
\item $\forall \alpha,\beta\in \omega_1\backslash \omega\,(\,f_\alpha|_{T_\alpha\cap T_\beta}=^*f_\beta|_{T_\alpha\cap T_\beta}\,).$
\end{enumerate}
 Given such family, we define $L^\xi_\alpha$ as $f^{-1}_\alpha[\{\xi\}]$ for all $\xi\in \omega_1$ and $\alpha> \xi$. Additionally, we let $\mathcal{L}^{\xi}=\langle L^\xi_\alpha\rangle_{\alpha>\xi}$
 \end{definition}

 Farah constructed a Haudorff coherent family of functions supported by an $\omega_1$-tower by forcing it and then appealing to Keisler's completeness Theorem for $L^\omega(Q)$. Unlike Talayco's result, the proof idea behind the construction of Farah does not work for any $\omega_1$-tower. This is due to the nature of Keisler's Theorem. 

As an application of the concept of AD representations, we  generalize the results of Talayco, Todor\v{c}evi\'c, Morgan and Farah regarding the existence of \say{big gap cohomology groups}. It is worth noting that there is no analogous to the $\aleph_0$ gap theorem for arbitrary $*$-lower semi-lattices. Namely, if there is a uniformizable AD system (also called strong-$Q$-sequence), then one can prove that there is a $*$-lower semilattice $\mathcal{T}$ isomorphic to $[\omega_1]^{<\omega}$ for which $\mathcal{G}(\mathcal{T})$ is a trivial group. In particular, we can not argue in the same way as Todor\v{c}evi\'c in order to conclude that the gap cohomology group has always cardinality $2^{\omega_1}$. Readers who wish to learn more about strong-$Q$-sequences are refered to \cite{strongqsequencesdavid} and \cite{Strongqjuris}.
\begin{theorem}\label{biggapcohomologytheorem} Let $(X,<,\wedge)$ be an $\omega_1$-like lower semi-lattice. Then there is a $*$-lower semi-lattice $\mathcal{T}$ isomorphic to $X$ with $|G(\mathcal{T})|=2^{\omega_1}.$
\begin{proof}Consider $X\times 2$ with the order given by $(x,i)<(y,j)$ if and only if $x<y$ and $i=j$. It is easy to see that $X\times 2$ is also $\omega_1$-like. In virtue of Theorem \ref{luzinreptheorem} we can take $\mathcal{T'}=\langle T'_{(x,i)}\rangle_{(x,i)\in X\times 2}$ and $\mathcal{A}=\langle A_{(x,i)}\rangle_{(x,i)\in X\times 2}$ so that $(\mathcal{T}',\mathcal{A})$ is an almost disjoint representation of $X\times 2$. For any $x\in X$ let $T_x=T'_{(x,0)}\cup T'_{(x,1)}. $ It is straightforward that the following properties hold for any $x,y \in X$:
\begin{enumerate}[label=$(\alph*)$]
\item $A_{(x,0)}\cup A_{(x,1)}\subseteq T_x$.

\item If $y\not\leq x$ then $A_{(y,0)}\cup A_{(y,1)}\subseteq^*T_y\backslash T_x$. In other words, $(A_{(y,0)}\cup A_{(y,1)})\cap T_x=^*\emptyset$.
\item $T_x\cap T_y=^*T_{x\wedge y}.$

\end{enumerate}
In particular, properties (b) and (c) imply that $\mathcal{T}=\langle T_x\rangle_{x\in X}$ is a $*$-lower semi-lattice isomorphic to $X$. We claim that $|G(\mathcal{T})|=2^{\omega_1}$. Trivially $|G(\mathcal{T})|\leq 2^{\omega_1}$, so we will only prove the other inequality. For this, let $S\in \mathscr{P}(X)$. . We will build recursively (using that $X$ is well-founded) a coherent subsystem $g_S:\mathcal{T}\longrightarrow \mathscr{P}(\omega)$ such that $A_{(x,0)}\subseteq g_S(T_x)$ for any $x\in X$ and:\\
\begin{center}\begin{minipage}{5cm} \begin{center} \textbf{(A)}\end{center} If $x\in S$, then $A_{(x,1)}\cap g_S(T_x)=^*\emptyset.$
\end{minipage}\hspace{2.5cm} \begin{minipage}{5cm}\begin{center} \textbf{(B)}\end{center} If $x\not\in S$, then $A_{(x,1)}\subseteq g_S(T_x)$.
\end{minipage}

\end{center}
Suppose that $y\in X$ and we have defined $g_S(T_x)$ for any $x<y$. By means of Lemma \ref{lemmacoherentsubsystempregap} we know that $\mathcal{G}(\mathcal{T}|_y,g_S|_{\mathcal{T}|_{y}})$ is pregap where $\mathcal{T}|_y=\langle T_x\rangle _{x<y}$. As there are no countable pregaps and $|(-\infty,y)|\leq \omega$, we conclude that there is $C\in \mathscr{P}(\omega)$ separating it.  Since $g_S(T_x)\subseteq T_x\subseteq^*T_y$ for any $x<y$, we may assume without loss of generality that $C\subseteq T_y$. Note that $C\cap T_x=^*g_S(T_x)$ for any $x<y$. We now define $g_S(T_y)$ by cases. If $y\in S$ define $g_S(T_y)$ as $(C\cup A_{(y,0)})\backslash A_{(y,1)}$. Otherwise define $g_S(T_y)$ as $C\cup A_{(y,0)}\cup A_{(y,1)}$. Since both $A_{(y,0)}$ and $A_{(y,1)}$ are almost disjoint with $T_x$ for all $x<y$ due to the point (b), it follows that $g_S(T_y)\cap T_x=^* C\cap T_x=^*g_S(T_x)$ for any such $x$. This finishes the recursion.

Now consider $\mathcal{S}\subseteq [X]^{\omega_1}$  a family of cardinality $2^{\omega_1}$ so that $S\Delta S'$ is uncountable for any two distinct $S,S'\in \mathcal{S}$. The key fact needed to finish the proof is that if $x\in S\Delta S'$ then $A_{(x,1)}\subseteq g_S\cdot g_{S'}(T_x)$ and $A_{(x,0)}\subseteq T_x\backslash (g_S\cdot g_{S'}(T_x))$.  In this way, any separation of the pregap  $\mathcal{G}(\mathcal{T},g_S\cdot g_{S'})$ would also separate the pregap $(A_{(x,1)},A_{(x,0)})_{x\in S\Delta S'}$ which is impossible since $\mathcal{A}$ is an inseparable family and $S\Delta S'$ is uncountable. This shows that $g_S$ and $g_{S'}$ are not cohomologous. Since $|S|=2^{\omega_1}$, then $|G(\mathcal{T})|=2^{\omega_1}$.
\end{proof}
\end{theorem}

\section{Finitizations of $\varclubsuit$ and CH from construction schemes}\label{sectionfinitizationcapturingaxioms}

Guessing principles are among the first things which come to mind when speaking about combinatorics over $\omega_1$. Some of the most utilized guessing principles are Jensen's $\Diamond$, the Continuum Hypothesis, and Ostszewski's $\varclubsuit$. These three axioms can be formulated in order to assert the existence of a sequence $\langle C_\gamma\rangle_{\gamma\in\text{Lim}}$ (where $\text{Lim}$ denotes the limit ordinals below $\omega_1$) of countable subsets of $\omega_1$ such that $C_\gamma\subseteq \gamma$ for each $\gamma$, and such that for any $S\in[\omega_1]^{\omega_1}$ and $\beta\in \omega_1$, the following property hold:
\begin{enumerate}[leftmargin=9em, itemsep=0.5em]
    \item[(Under $\Diamond$)] $\{\gamma\in \text{Lim}\,:\,\gamma\cap S=C_\gamma\text{ and }\sup(C_\gamma)=\gamma\}$ is stationary. A sequence satisfying this property its called a \emph{$\Diamond$-sequence}.
    \item[(Under $\varclubsuit$)] $\{\gamma\in \text{Lim}\,:\,C_\gamma\subseteq S\text{ and }\sup(C_\gamma)=\gamma\}$ is stationary. A sequence satisfying this property its called a \emph{$\varclubsuit$-sequence}.
    \item[(Under CH)] $\{\gamma\in \text{Lim}\,:\,C_\gamma\sqsubseteq S\text{ and }\beta\leq \sup(C_\gamma)\}$ is stationary.
\end{enumerate}
It can be seen at a first glance that these principles are related. For example, given a sequence $\langle C_\gamma\rangle_{\gamma\in \Lim}$, the intersection of the two latter sets described above is equal the former one when we take $\beta=0$. From this fact, it already follows that $\Diamond$ implies both CH and $\varclubsuit$. It is not completely direct that $\varclubsuit+CH$ imply $\Diamond$, even if there is one single sequence testifying these two axioms.  Nevertheless, such implication indeed hold, and in some sense, this tells us that the interplay between these axioms goes as far as possible.\\
In \cite{schemenonseparablestructures}, Todor\v{c}evi\'c introduced a generalization of Jensen's $(\omega,1)$-gap morasses which he called construction schemes. These objects are special families of finite subsets of $\omega_1$ which serve as a tool for building uncountable objects by means of finite approximations. In that same paper, he introduced the capturing axiom $FCA$. Such axiom states the existence of construction schemes with certain guessing properties which can be interpreted as finite dimensional versions of  the properties possessed by a $\Diamond$-sequence. In fact, it turns out that $FCA(part)$ is implied by $\Diamond$. As it has been proved in papers such as \cite{schemescruz}, \cite{irredundantsetsoperator}, \cite{banachspacescheme}, \cite{lopezschemethesis}, \cite{treesgapsscheme}  and \cite{schemenonseparablestructures}, this capturing axiom is sufficient to imply the existence of many of the most well-known and studied objects whose existence already follows from Jensen's principle. There are variations of the axiom $FCA(part)$ called $CA_n(part)$(for each natural number $2\leq n$) and $CA(part)$  as well as their non \say{part} variants $FCA$, $CA_n$ and $CA$ whose hierarchy is shown in the following diagram.\\
\begin{center}\small
\begin{tikzcd}
CA_2(part) \arrow[d, Rightarrow] &CA_3(part) \arrow[l, Rightarrow] \arrow[d, Rightarrow]&\dots \arrow[l, Rightarrow] \arrow[d, Rightarrow] & CA(part)\arrow[l, Rightarrow]\arrow[d, Rightarrow] & FCA(part)\arrow[d, Rightarrow] \arrow[l, Rightarrow] & \Diamond \arrow[l,Rightarrow] \\
CA_2 & CA_3 \arrow[l, Rightarrow] &\dots\arrow[l, Rightarrow] & CA\arrow[l, Rightarrow] & FCA\arrow[l, Rightarrow]
\end{tikzcd}
\end{center}

The arrows in this diagram are not reversible. Hence, this collection of axioms seems to be a good framework for studying the strength of the consequences of $\Diamond$. Of course, there are already other kinds of axioms which serve a similar purpose.  For example, the parametrized $\Diamond$-principles introduced in \cite{ParametrizedDiamonds} by D\v{z}amonja, Hru\v{s}\'ak and Moore. \\
The main goal of this section is to describe finite dimensional versions of $\varclubsuit$ and CH, called $CA^\rho$ and $CA^\Delta$ (and their respective variants), which play an analogous roll in relation to the capturing axioms for construction schemes as the one played by their \say{infinite dimensional} counterparts with respect to the $\Diamond$-principle.\\
We assume that the general audience is yet not familiar with the theory of construction schemes. For that reason, we will dedicate the first part of this section to recall some of the basic definitions and results regarding this theory. After doing so, we will be in the position to state the axioms  $CA^\rho$ and $CA^\Delta$ in a precise way.\\
 A \emph{type} is a sequence $\tau=\langle m_k,n_{k+1},r_{k+1}\rangle_{k\in\omega}$ of triplets of natural numbers such that the following conditions hold for any $k\in\omega$:\\
\begin{enumerate}[label=$(\alph*)$,itemsep=0.5em]
\begin{minipage}{5cm}
\item $m_0=1,$
\item $n_{k+1}\geq 2,$
\end{minipage}
\begin{minipage}{5cm}
\item $m_k>r_{k+1}, $
\item $m_{k+1}=r_{k+1}+(m_k-r_{k+1})n_{k+1}.$
\end{minipage}\\
\end{enumerate}
We say that type is \textit{good} if for each $r\in\omega$, there are infinitely many $k$'s for which $r=r_k$. Throughout this text, we will only work with good types. We say that  partition of $\omega$, namely $\mathcal{P}$, is \emph{compatible with $\tau$} if $\tau$ is good when restricted to each member of $\mathcal{P}$. That is, for each $P\in \mathcal{P}$ and every $r\in \omega$ there are infinitely many $k\in P$ for which $r_k=r$.

Given a set $X$ of ordinals, a construction scheme (of type $\tau$) over $X$ is a family $\mathcal{F}\subseteq \text{Fin}(X)$ which: Is cofinal in $(\text{Fin}(X),\subseteq )$, any member of $\mathcal{F}$ has cardinality $m_k$ for some $k\in\omega$, and furthermore, if we put $\mathcal{F}_k:=\{ F\in \mathcal{F}\,:\,|F|=m_k\}$, then the following two properties are satisfied for each $k\in\omega$:\vspace{0.5em}
\begin{enumerate}[label=(\roman*),itemsep=0.5em]
\item $\forall F,E\in \mathcal{F}_k\big(\; E\cap F\sqsubseteq E,F\;\big)$,
\item $\forall F\in \mathcal{F}_{k+1}\;\exists F_0,\dots,F_{n_{k+1}-1}\in \mathcal{F}_k$ such that $$F=\bigcup\limits_{i<n_{k+1}}F_i.$$
Moreover, $\langle F_i\rangle_{i<n_{k+1}}$ forms a $\Delta$-system with root $R(F)$ such that  $|R(F)|=r_{k+1}$ and $R(F)<F_0\backslash R(F)<\dots < F_{n_{k+1}-1}\backslash R(F).$
\end{enumerate}\vspace{0.5em}
Note that for a given $F\in \mathcal{F}_{k+1}$, each of the $F_i$'s mentioned above can be written as $F[r_{k+1}]\cup F[\,[a_i,a_i+1)\,]$ where $a_i=r_{k+1}+i\cdot(m_{k+1}-m_k)$. In particular, this is saying that the family $\langle F_i\rangle_{i<n_{k+1}}$ is unique, so  we call it the \emph{canonical decomposition} of $F$. If $\mathcal{F}$ is a construction scheme over $X$, we  call $X$ \emph{the domain} of $\mathcal{F}$ and we write it as $dom(\mathcal{F})$. In practice, we are only interested in construction schemes over $\omega_1$. However, in order to build such schemes it is convenient to work also with schemes with both finite and countable domains. We will frequently use the following fact. 
\begin{proposition}\label{propuniquescheme}Let $\tau$ be a good type and $X\in \Fin(\omega_1)\cup \{\omega\}$. There is a construction scheme (of type $\tau$) over $X$ if and only if either $X=\omega$ or $|X|=m_k$ for some $k\in\omega$. In such case the scheme is unique, so we will call it $\mathcal{F}(X)$.
\end{proposition}
\begin{rem}Suppose that we have a construction scheme $\mathcal{F}$ and $F\in \mathcal{F}$. According to the Proposition \ref{propuniquescheme}, there is a unique construction scheme over $F$, namely $\mathcal{F}(F)$. It is easy to see that in this case, $\mathcal{F}(F)=\{G\in \mathcal{F}\,:\,G\subseteq F\}$. Particularly, this means that $\mathcal{F}=\bigcup\limits_{F\in \mathcal{F}} \mathcal{F}(F)$.
\end{rem}

Each construction scheme $\mathcal{F}$ is associated with some natural functions of countable codomain. The first of such functions is the one derived from the cofinality condition in the definition of a scheme. Namely, given $\alpha,\beta\in dom(\mathcal{F})$ there is $F\in \mathcal{F}$ such that $\{\alpha,\beta\}\subseteq F$. Hence, we can define $\rho:dom(\mathcal{F})^2\longrightarrow \omega$ as:$$\rho(\alpha,\beta)=\min (k\in\omega\,:\,\exists F\in \mathcal{F}\,(\{\alpha,\beta\}\subseteq F )\,).$$
For each finite $A\in dom(\mathcal{F})$, we also define $$\rho^A=\max(\rho[A^2])=\max(\rho(\alpha,\beta)\,:\,\alpha,\beta\in A).$$
It is not hard to see that $\rho^F=n$ for each $F\in \mathcal{F}_n$.
The most important feature of the function $\rho$ is that it is an ordinal metric. This means that it satisfies the properties stated in the following lemma.
\begin{lemma}Let $\mathcal{F}$ be a construction scheme. The following properties hold for any $\alpha,\beta,\gamma\in dom(\mathcal{F})$ and each $k\in\omega$:\vspace{0.5em}
\begin{enumerate}[label=$(om_{\arabic*})$, itemsep=0.5em]
\item $\rho(\alpha,\beta)=0$ if and only if $\alpha=\beta$.
\item $\rho(\alpha,\beta)=\rho(\beta,\alpha).$ 
\item If $\alpha\leq \min(\beta,\gamma)$, then $\rho(\alpha,\beta)\leq \max(\,\rho(\alpha,\gamma),\rho(\beta,\gamma)\,)$.
\item $\{\xi\leq\alpha\,:\,\rho(\alpha,\xi)\leq k\}$ is finite.
\end{enumerate}
\end{lemma}
Given $\alpha\in dom(\mathcal{F})$ and $k\in\omega$ we can define the $k$-closure of $\alpha$ as $(\alpha)_k:=\{\xi\leq \alpha\,:\rho(\alpha,\xi)\leq k\}$ and $(\alpha)_k^-:=(\alpha)_k\backslash \{\alpha\}$. Note that property $(om_4)$ is saying that all the $k$-closures are finite. It is a useful fact that for any $k\in\omega$ (with $m_k\leq |dom(\mathcal{F})|$) there is at least one $F\in \mathcal{F}_k$ such that $\alpha\in F$. Even more, for any such $F$ we have the equalities: $$F\cap (\alpha+1)=(\alpha)_k$$
$$F\cap \alpha=(\alpha)^-_k.$$
Given $\alpha\in dom(\mathcal{F})$, the $k$-cardinality function $\lVert \alpha\rVert_{\_}:\omega\longrightarrow \omega$ is defined as:
$$\lVert \alpha\rVert_k=|(\alpha)^-_k|.$$
It is helpful to remember that if $F\in \mathcal{F}_k$ is such that $\alpha\in F$, then $F(\lVert \alpha \rVert_k)=\alpha$. Recall that if $f,g:\omega\longrightarrow\omega$ are distinct,  we can define $\Delta(f,g):=\min (\,k\in \omega\,:\,f(k)\not=g(k)\,)$. We put $\Delta(f,g):=\omega$ whenever $f=g$. In the case of  construction schemes, this leads to the definition of the function $\Delta:dom(\mathcal{F})^2\longrightarrow \omega+1$ which is defined as: $$\Delta(\alpha,\beta)=\Delta(\lVert \alpha\rVert_{\_},\lVert \beta\rVert_{\_}).$$
In the following lemma, we describe the most basic properties of this function as well as its relation with the $\rho$-function.
\begin{lemma}Let $\mathcal{F}$ be a construction scheme and let $\alpha,\beta,\delta\in dom(\mathcal{F})$. Then the following properties hold:\vspace{0.5em}
\begin{enumerate}[label=$(dp_{\arabic*})$, itemsep=0.5em]
    \item If $\alpha<\beta$, then $\lVert \alpha\rVert_k<\lVert \beta\rVert_k$ for each $k\geq \rho(\alpha,\beta)$. In particular, $\Delta(\alpha,\beta)\leq \rho(\alpha,\beta)$ whenever $\alpha\not=\beta$.
    \item If $\Delta(\alpha,\beta)<\Delta(\beta,\delta)$, then $\Delta(\alpha,\beta)=\Delta(\alpha,\delta)$.
\end{enumerate}
\end{lemma}
The last canonical function that we will present is the one which is related to the property $(ii)$ in the definition of a scheme. Namely, given $\alpha\in dom(\mathcal{F})$ we define the function $\Xi_\alpha:\omega\longrightarrow \{-1\}\cup \omega$ as follows;  If $1\leq k\in\omega$ and $m_k\leq |X|$, then there is $F\in \mathcal{F}_k$ so that $\alpha\in F$. According to the property $(ii)$ in the definition of a construction scheme, we have that either $\alpha\in R(F)$ or there is a unique $i<n_k$ such that $\alpha\in F_i\backslash R(F)$. We then define $$\Xi_\alpha(k):=\begin{cases}-1&\text{ if }\alpha\in R(F)\\
i&\text{ if }\alpha\in F_i\backslash R(F)\end{cases}$$ It can be proved that this definition does not depend on the choice of $F.$ Now, if either $k=0$ or $m_k>|X|$, we simply define $\Xi_\alpha(k)$ as $0$.\\
The next lemma relates the functions $\rho$, $\Xi$ and $\Delta$ in a very nice way. In fact, this is, by experience, the most utilized result regarding the basic theory of construction schemes.

\begin{lemma}\label{lemmaxi}Let $\mathcal{F}$ be a construction scheme, $\alpha<\beta\in dom(\mathcal{F})$ and $1\leq k\in \omega$ be such that $m_k\leq |dom(\mathcal{F})|$. Then:\vspace{0.5em}
\begin{enumerate}[label=$(xi_{\alph*})$, itemsep=0.5em]
\item If $k<\Delta(\alpha,\beta)$, then  $\Xi_\alpha(k)=\Xi_\beta(k).$
\item If $k=\rho(\alpha,\beta)$, then $0\leq \Xi_\alpha(k)<\Xi_\beta(k).$
\item If $k>\rho(\alpha,\beta)$, then either $\Xi_\alpha(k)=-1$ or $\Xi_\alpha(k)=\Xi_\beta(k).$ 
\item If $k=\Delta(\alpha,\beta)$ then $0\leq \Xi_\alpha(k)\not=\Xi_\beta(k)\geq 0.$ 
\end{enumerate}
\end{lemma}
As suggested in the previous lemma, there is not much to be said about the functions $\Xi_\alpha$ and $\Xi_\beta$  for the $k$'s which live in the open interval $(\Delta(\alpha,\beta),\rho(\alpha,\beta))$. On the other hand, the relation between $\Xi_\alpha$ and $\Xi_\beta$ can be fully described for $k$'s outside the afforded mentioned interval. At this point is where the notion of capturing comes into play. Naively speaking, a construction scheme over $\omega_1$ is capturing (up to some degree) if for any uncountable subset $S$ of $\omega_1$ there are distinct $\alpha,\beta$ for which the gap between $\Delta(\alpha,\beta)$ and $\rho(\alpha,\beta)$ is non-existent. The importance of this property is that, in this case, Lemma \ref{lemmaxi} can be used at its full power.

In order to motivate the relation between the already studied capturing axioms and the ones which we study in this manuscript, we will present them in a way which reassembles the presentation of $\Diamond$, $\varclubsuit$ and CH given at the beginning of this section.
\begin{definition}[$\ast$-Captured $\Delta$-systems]\label{capturedsystemdef} Consider $\mathcal{F}$ be a construction scheme, $X=dom(\mathcal{F})$, and $\ast$ a family of functions from $X^2$ to $\omega$. Let $1\leq n,m\in\omega$ and $\mathcal{D}=\langle D_i\rangle _{i<n}\subseteq [X]^{m}$ be a root-tail-tail $\Delta$-system with a root $R$ of cardinality $r$. Given $l\in \omega$, we will say that $\mathcal{D}$ is $\ast$-\emph{captured} at level $l$ if \\ 
\begin{center}
\begin{minipage}{10cm}\begin{center}\textbf{(I)}\\
For all $i<n$ and each $a<m$, $$\Xi_{D_i(a)}(l)=\begin{cases}-1&\text{if }a\leq r\\
i&\text{if }a>r
\end{cases}$$
\end{center}
\end{minipage}
\end{center}
and for every $\nu\in \ast$ the following condition hold:
\begin{center}\begin{minipage}{10cm}
\begin{center} \textbf{(II$_\nu$)} \\
 For all $i<j<n$ and $r\leq a<m$,
$$\nu(D_i(a),D_j(a))=l.$$\end{center}
\end{minipage}\vspace{0.5cm}
\end{center}

 We denote the family of all $\ast$-captured $\Delta$-systems  at level $l$ as $\text{Cap}^\ast_l(\mathcal{F})$.
Finally, if $\mathcal{D}$ is $\ast$-captured at level $l$ and $n=n_l$, we will say that $\mathcal{D}$ is $\ast$-\emph{fully captured}.  The family of all $\Delta$-systems which are $\ast$-fully captured at level $l$ will be denoted as $\text{FCap}^\ast_l(\mathcal{F})$. 
\end{definition}

Suppose that $\mathcal{F}$ and $\ast$ are as in the previous definition and that $\omega_1=dom(\mathcal{F})$. We say that $\mathcal{F}$ is 
\emph{$(\mathcal{P},n)$-$\ast$-capturing} (for some $n\in\omega$ and a partition $\mathcal{P}$  of $\omega$) if for any uncountable $\mathcal{S}\subseteq \text{Fin}(\omega_1)$ and each $P\in\mathcal{P}$, the set $$\{l\in P\,:\,[\mathcal{S}]^n\cap \text{Cap}^\ast_l(\mathcal{F})\not=\emptyset\}$$ is infinite. If $\mathcal{F}$ is $(\mathcal{P},n)$-$\ast$-capturing for any $2\leq n \in\omega$, we will say that it is \emph{$\mathcal{P}$-capturing}. Finally, we will say that $\mathcal{F}$ is fully \emph{$\mathcal{P}$-$\ast$-capturing} if for any uncountable $\mathcal{S}\subseteq \text{Fin}(\omega_1)$ and $P\in\mathcal{P}$, the set $\{ l\in P\,:\,\text{Fin}(S)\cap\text{FCap}^\ast_l(\mathcal{F})\not=\emptyset\}$ is infinite. Whenever $\mathcal{P}$ is the trivial partition $\{\omega\}$, we will omit it from the notation.
\begin{rem}
It is worth to highlight how the $\ast$-capturing schemes reassemble the definition of a $\Diamond$-sequence. Namely, for $\ast$-fully capturing construction schemes the sequence $\langle C_\gamma\rangle_{\gamma\in \text{Lim}}$ is changed by $\langle \text{FCap}^\ast_l(\mathcal{F})\rangle_{l\in\omega}$, the objects that we are trying to guess are uncountable $\mathcal{S}\subseteq \text{Fin}(\omega_1)$ instead of uncountable $S\in[\omega_1]^{\omega_1}$, the guessing condition  $\gamma\cap S=C_\gamma$  turns into $\text{FCap}^\ast_l(\mathcal{F})\cap \text{Fin}(\mathcal{S})\not=\emptyset$, and the word stationary is substituted by the word infinite.
\end{rem}
By now, we have presented all the necessary information needed to state the axioms which gave rise to this manuscript. In order to avoid unnecessary verbatim, let us first fix $\ast\in\{\,\{\Delta\},\{\rho\},\{\Delta,\rho\}\,\}$:\\\\
\noindent
{\bf Fully $*$-Capturing Axiom} {[\bf FCA$^\ast$]}: There is a fully  $\ast$-capturing construction scheme over $\omega_1$ of every possible good type.\\\\
{ \bf Fully Capturing Axiom with Partitions [FCA$^\ast$(part)]}: There is a fully $\mathcal{P}$-$\ast$-capturing construction scheme over $\omega_1$ for every good type $\tau$  and each partition $\mathcal{P}$ compatible with $\tau$.\\\\
{\bf $n$-$\ast$-Capturing Axiom [CA$^\ast_n$]}: For any $n'\leq n$, there is an $n'$-$\ast$-capturing construction scheme over $\omega_1$ of every possible good type satisfying that $n'\leq n_k$ for each $k\in \omega\backslash1$.\\\\
{\bf $n$-$\ast$-Capturing Axiom with Partitions [CA$^\ast_n$(part)]}: For any $n'\leq n$, there is a $(\mathcal{P},n')$-$\ast$-capturing construction scheme over $\omega_1$ for every good type $\tau$ satisfying that $n'\leq n_k$ for each $k\in \omega\backslash 1$  and each partition $\mathcal{P}$ compatible with $\tau$.\\\\
{\bf $\ast$-Capturing Axiom [CA$^\ast$]}: CA$^\ast_n$ holds for each $n\in\omega$ and there is a $\ast$-capturing construction scheme over $\omega_1$ for every good type satisfying that the sequence $\langle n_{k+1}\rangle_{k\in\omega}$ is non-decreasing and unbounded.\\\\
{\bf Capturing Axiom with partitions [CA$^\ast$(part)]}: CA$^\ast_n$(part) holds for each $n\in\omega$ and there is a $
\mathcal{P}$-$\ast$-capturing construction scheme over $\omega_1$ for every good type $\tau$ satisfying that the sequence $\langle n_{k+1}\rangle_{k\in\omega}$ is non-decreasing and unbounded and each partition $\mathcal{P}$ compatible with $\tau$.\\
\begin{remark} The axioms presented above coincide with the \say{capturing axioms} already introduced in \cite{schemenonseparablestructures}. provided that   $*=\{\rho,\Delta\}$.  In this particular case, we suppress the $*$ from the notation. 
\end{remark}

The two following theorems will be proved in sections \ref{sectionclub} and \ref{sectionch} respectively. It is important to remark that both $\text{CA}^\Delta_2$ and $\text{CA}^\rho_2$ hold without assuming any extra axioms. 

\begin{theorem}[Under CH]\label{fcadeltach} $FCA^\Delta$ holds.
\end{theorem}
\begin{theorem}[Under $\varclubsuit$ principle]\label{carhoclub} $CA^\rho$ holds.
\end{theorem}
It will be clear that simple modifications of the proofs of the theorems above, actually yield that $FCA^\Delta(part)$ and $CA^\rho(part)$   follow from CH and $\varclubsuit$ respectively.

\section{Proof of Theorem \ref{omegainseptheorem}}\label{sectionproofomegainseptheorem}
In this section we will prove Theorem \ref{omegainseptheorem}, namely, that if either $\varclubsuit$ or CH hold, then there is an $\omega$-inseparable $AD$ family which codes every $\omega_1$-like order. What we will actually do is to prove via the theorems \ref{luzinjonestheorem} and \ref{luzinreptheorem} that if there is either a $n$-$\Delta$-capturing or $n$-$\rho$-capturing construction scheme, then there is an $n$-inseparable family which codes any $\omega_1$-like order (here, $n\in \omega+1$). 

\begin{theorem}\label{luzinjonestheorem}Let $\mathcal{F}$ be a scheme of type $\langle m_k,n_k,r_k\rangle_{k\in \omega}$ and $n\in \omega$. Given $k\in \omega$, consider $N_k= m_{k-1}^{n_k}\times \{k\}$ and $N=\bigcup\limits_{k\in \omega} N_k$. Given $\alpha\in \omega_1$ and $k\in \omega$ we define $$A^k_\alpha:=\{\sigma\in N_k\,:\,\Xi_\alpha(k)\geq 0\text{ and }\sigma(\Xi_\alpha(k))=\lVert \alpha\rVert_{k-1}\}$$
and $A_\alpha:=\bigcup\limits_{k\in \omega} A_k$.  Then $\mathcal{A}_\mathcal{F}=\langle A_\alpha\rangle_{\alpha\in \omega_1}$ is an inseparable $AD$-family. Furthermore, if  either $\mathcal{F}$ is $n$-$\rho$-capturing or $n$-$\Delta$-capturing, then $\mathcal{A}_\mathcal{F}$ is an $n$-inseparable. In particular, if $\mathcal{F}$ is either $\Delta$-capturing or $\rho$-capturing, then $\mathcal{A}_\mathcal{F}$ is $\omega$-inseparable.
\begin{proof}Each $A_\alpha$ is infinite because $\Xi_\alpha(k)\geq 0$ for infinitely many $k$'s. It is also easy to see that for distinct $\alpha$ and $\beta$ we have that$$
    A_\alpha\cap A_\beta\subseteq \bigcup\limits_{k\leq \rho(\alpha,\beta)} (A_\alpha^k\cap A_\beta^k).
$$
Therefore, $A_\mathcal{F}$ is an AD family. Before proving the second part note that if $D\in[\omega_1]^n$ and $k\in \omega$ are such that $\Xi_{D(i)}(k)=i$ for each $i$, then $$
   \{\sigma\in N_k\,:\,\forall i<n\,( \sigma(i)=\lVert D(i)\rVert_{k})\}\subseteq \bigcap\limits_{i<n} A^k_{D(i)}.
$$
So in particular, $\big(\bigcap\limits_{i<n} A^k_{D(i)}\big)\cap N_k\not=\emptyset$.\\
Assume that $\mathcal{F}$ is  $n$-$\rho$-capturing (resp. $n$-$\rho$-capturing) and let $\langle I_i\rangle_{i<n}$ be a sequence of uncountable disjoint subsets of $\omega_1$. We shall prove that $( \langle A_\alpha\rangle_{\alpha\in I_i}\,:\,i<n)$ is an $n$-gap. Assume towards a contradiction that 
this is not the case and let $c:n\longrightarrow \mathscr{P}(N)$ be a witness of this fact. By refining each $I_i$, we may assume without loss of generality that for each $i<n$ there is $k_i\in \omega$ so that $A_\alpha^k\subseteq c(i)$ for all $\alpha\in I_i$. As $\bigcap\limits_{i<n} c(i)$ is finite, there is $k>\max(k_i\,:i<n)$ for which this intersection is disjoint from $A_l$ for each $l\geq k$. Now, we construct a family $\mathcal{D}\subseteq[\omega_1]^n$ be an uncountable family of pairwise disjoint sets so that $D(i)\in I_i$ for all $D\in \mathcal{D}$ and $i<n$. Siince $\mathcal{F}$ is $n$-$\rho$ -capturing (resp. $n$-$\Delta$-capturing) there exists a subfamily $\langle D_i\rangle_{i<n}$ of $\mathcal{D}$ and $l>k$ for which $\langle D_i\rangle_{i<n}$ is $\rho$-captured (resp.$\Delta$-captured) at level $l$. Since the $D_i$'s are pairwise disjoint, we have that $\Xi_{D_i(i)}(l)=i$ for each  $i<n$. Therefore, $\emptyset\not=\big(\bigcap\limits_{i<n}A^k_{D_i(i)} \big)\cap N_l\subseteq \big(\bigcap\limits_{i<n} c(i)\big)\cap N_l$. This is a contradiction so we are done.
\end{proof}
    
\end{theorem}

Now, we will prove that $\mathcal{A}_\mathcal{F}$ codes any $\omega_1$-like order. For that, we need the following easy lemma.
\begin{lemma}\label{bijectionlemma}Let $(X,\leq)$ be $\omega_1$-like. Then there is a bijection $\phi:X\longrightarrow \omega_1$ with $\phi(x)<\phi(y)$ for all $x<y\in X.$
\end{lemma}

Let $(X,\leq)$ $\omega_1$-like and $\phi:X\longrightarrow \omega_1$ as above. In the following lemma, in some sense, we will pull back a construction scheme from $\omega_1$ to some sort of \say{construction scheme} on $X$. This suggests that the theory of construction schemes could be generalized to other partial orders of size $\omega_1$. The extend and utility of this idea remains mainly unexplored.

\begin{lemma}\label{Mklemma}Let $\mathcal{F}$ be a scheme and let $X$ and $\phi$ be as in Lemma \ref{bijectionlemma}. For each $x\in X$ and $k\in\omega$ define $M^k_x=\{\, z\leq x\,:\, \phi(z)\in (\phi(x))_k\,\}$. The following properties hold for each $x,y\in X$ and $k\in\omega$:
\begin{enumerate}[label=$(\arabic*)$]
\item  If $\inf(x,y)$ exists and $k>\rho^{\phi[\{x,y,\inf(x,y)\}]}$, then $M^k_x\cap M^k_y=M^k_{\inf(x,y)}.$
\item[$(\frac{3}{2})$] If $y\leq x$ and $k>\rho^{\phi[\{x,y\}]}$, then $M^k_y\subseteq M^k_x$.
\item If $y\not\leq x$, then $y\in M^k_y\backslash M^k_x.$
\item If $x$ is succesor-like and $k>\rho^{\phi[pred(x)\cup\{x\}]}$, then $M^k_x\backslash\big( \bigcup\limits_{z\in pred(x)} M^k_z \big)=\{x\}$.
\item If $x$ and $y$ are incompatible,  then $M^k_x\cap M^k_y=\emptyset.$
\end{enumerate}
\begin{proof}
The points (2) and (4) are trivial and the point ($\frac{3}{2}$) follows directly from (1). Therefore,  we will only prove (1) and (3).

\begin{claimproof}[Proof of  $(1)$] Since $\phi$ preserves the order then $\phi(\inf(x,y))\leq \min(\phi(x),\phi(y))$ which implies that $(\phi(\inf(x,y))_k\sqsubseteq (\phi(x))_k\cap (\phi(y))_k$. In this way, if $z\in M^k_{\inf(x,y)}$ then $z\leq x,y$ by definition of the infimum, and $z$ belongs to both $(\phi(x))_k$ and $(\phi(y))_k$ because $z\in (\phi(\inf(x,y)))_k$. Therefore $z\in M^k_x\cap M^k_y$. On the other hand, if $z\in M^k_x\cap M^k_y$ then $z\leq x$ and $z\leq y$, so $z\leq \inf(x,y)$. Furthermore, $\phi(z)\in (\phi(x))_k\cap (\phi(y))_k$ and $\phi(z)\leq \phi(\inf(x,y))$. As $(\phi(\inf(x,y)))_k$ is an initial segment of such intersection then $\phi(z)$ belongs to it. Consequently $z\in M^k_{\inf(x,y)}$. 
\end{claimproof}
\begin{claimproof}[Proof of $(3)$] Note that since $x\not\leq z$ for any $z\in pred(x)$ then $x\in M^k_x\backslash M^k_z$ for any such $z$. This proves the inclusion from right to left. To show that the one from left to right also holds, let $$w\in M^k_x\backslash\big( \bigcup\limits_{z\in pred(x)} M^k_z \big).$$
Suppose towards a contradiction that $w<x$ and consider $z\in pred(x)$ with $w\leq z$. Since $k>\rho(\phi(x),\phi(z))$ we conclude that $(\phi(z))_k\sqsubseteq (\phi(x))_k$. But then  $\phi(w)\in (\phi(z))_k$ because $\phi(w)\leq \phi(z)$ and $\phi(w)\in (\phi(x))_k$. Consequently $w\in M^k_z$ which is a contradiction. This finishes the proof.
  \end{claimproof}  

\end{proof}
\end{lemma}

\begin{theorem}\label{luzinreptheorem}Let $\mathcal{A}_\mathcal{F}$ be the $AD$ family constructed in Theorem \ref{luzinjonestheorem} and let $(X,\leq)$ an $\omega_1$-like order. Then there is a family $\mathcal{T}=\langle T_x\rangle_{x\in X}\subseteq [\omega]^{\omega}$ and a re-indexing of $\mathcal{A}$ as $\langle\hat{A}_x\rangle_{x\in X}$ so that $(\mathcal{T},\mathcal{A})$ is an almost disjoint representation of $X$.
\begin{proof}Let $\phi:X\longrightarrow \omega_1$ be as in Lemma \ref{bijectionlemma} and
$ M^k_x$ be as in Lemma 
\ref{Mklemma} for any $x\in X$ and $k\in \omega$. Fix $x\in X$ and let $\hat{A}_x=A_{\phi(x)}.$
This defines the re-indexing of $\mathcal{A}$. Now, for each $k\in\omega$ and $x\in X$, let 
$$T^k_x=\bigcup \{\,A_{\phi(z)}^{k+1}\,:\,z\in M^k_x\,\}=\bigcup\{\,A^{k+1}_\xi\,:\,\xi\in \phi[M^k_x]\,\}$$
$$T_x=\bigcup\limits_{k\in\omega}T^k_x.$$
We claim that $(\mathcal{T},\mathcal{A})$ is an almost disjoint representation of $X$. In the following paragraphs we will show that the points $(a)$, $(b)$, $(c)$, $(d)$ and $(e)$ of Definition \ref{Luzinrepdef} are satisfied for such pair. 
\begin{claimproof}[Proof of $(a)$]Given $k>0$ we have that $A^k_{\phi(x)}\subseteq T^{k-1}_x$ because $x$ is always an element of $M^k_x$. Therefore $\hat{A}_x=A_{\phi(x)}\subseteq T_x$.   
\end{claimproof}
 Before proving the remaining points let us fix some notation. Given $P\in[\omega_1]^{<\omega}$ and $k\in\omega$, let $A^k_P=\bigcup\{A^k_\alpha\,:\,\alpha\in P\}$. Note that if $P,Q\in [\omega_1]^{<\omega}$ and $k>\rho^{P\cup Q}$, then we have the following:

 $$A^k_P\cap A^k_Q=A^k_{P\cap Q},$$
$$A^k_P\cup A^k_Q=A^k_{P\cup Q},$$
$$A^k_{P}\backslash A^k_Q=A^k_{P\backslash Q}.$$

We will apply these three equalities for finite sets of the form $\phi[M^k_x]$. More precisely, observe that if $x,y\in X$ and $k>\rho(\phi(x),\phi(y))$ then $\rho^{\phi[M^k_x]\cup \phi[M^k_y]}<k+1$. This is because if we take  $F\in \mathcal{F}_k$ with $\phi(x),\phi(y)\in F$, then $\phi[M^k_x]\cup \phi[M^k_y]=\phi[M^k_x\cup M^k_y]\subseteq (\phi(x))_k\cup (\phi(y))_k\subseteq F.$ 

\begin{claimproof}[Proof of $(c)$]Let $x,y\in X$ and suppose that $\inf(x,y)$ exists.  Take an arbitrary $k>\rho^{\phi[\{x,y,\inf(x,y)\}]}$. By means of the point $(1)$ of Lemma  \ref{Mklemma}, $M^k_x\cap M^k_y=M^k_{\inf(x,y)}$. In this way: \begin{align*}
T^k_x\cap T^k_y =A^{k+1}_{M^k_x}\cap A^{k+1}_{M^k_y}
=A^{k+1}_{M^k_x\cap M^k_y}
=A^{k+1}_{M^k_{\inf(x,y)}}=T^k_{\inf(x,y)}.
\end{align*}
Therefore, $T_x\cap T_y=\bigcup\limits_{k\in\omega} \big( T^k_x\cap T^k_y\big)=^*\bigcup \{T^k_{\inf(x,y)}\,:\, k>\rho^{\phi[x,y,\inf(x,y)]}\}=^*T_{\inf(x,y)}.$
\end{claimproof}
The remaining parts of the theorem are proved in a completely similar way. Because of this, we leave the calculations to the reader. 
\end{proof}
\end{theorem}

   \section{Entangled sets and $\sigma$-monotone spaces}\label{sectionentangled}
One of Cantor's most famous theorems is that any two countable dense total orders without endpoint are isomorphic. Particularly, this result can be applied to any two countable dense subsets of the reals. Thus, it is natural to ask whether Cantor's theorem can be extended to higher infinities (inside $\mathbb{R}$). 
\begin{definition}[$\kappa$-dense sets] Let $\kappa$ be an infinite cardinal. We say that $D\subseteq \mathbb{R}$ is \textit{$\kappa$-dense} if $|(a,b)\cap D|=\kappa$ for all $a<b\in D$.
    
\end{definition}
In \cite{baumgartneraxiom}, Baumgartner proved that it is consistent that any two $\omega_1$-dense sets of reals are isomorphic. This assertion is now known as the Baumgartner's axiom for $\omega_1$. More generally, if $\kappa$ is an infinite cardinal then we can define:\\

\noindent
{\bf Baumgartner's Axiom for $\kappa$ [$BA(\kappa)$]}: Any two $\kappa$-dense sets of reals are isomorphic.\\

\noindent
One of the ways of proving Cantor's theorem is by defining the forcing of finite approximations of isomorphisms and making use of the Sikorski's Lemma. Hence, it is natural to think that maybe $BA(\kappa)$ has some relation with  $MA$. However, this is not the case. In \cite{MAdoesnotImplyBA}, Abraham and Shelah showed that $MA$ does not imply $BA(\omega_1)$. For that, they introduced the objects which we  will study in this section. Namely, the entangled sets. Readers interested in learning more about them  may also look at \cite{WhyYcc}, \cite{guzmanentangled}, \cite{AsperoMotaEntangled}, \cite{RemarksonChainConditionsinProducts}, \cite{PartitionProblems} and \cite{EntangledCohen}.    
\begin{definition}[Realization]Let $k\in\omega$ and $t:k\longrightarrow \{<,>\}$. Given a total order $(X,<)$ and $a,b\in [X]^k$, we say that $(a,b)$ \textit{realizes} $t$ if $$a(i)\:t(i)\:b(i)$$
for each $i<k.$ By $T(a,b)$ we denote the unique $t$ which is realized by $(a,b).$
\end{definition}  

\begin{definition}[Entangled set]Let $(X,<)$ be a partial order, $\mathcal{E}\in [X]^{\omega_1}$ and $k\in\omega.$ We say that:\begin{itemize}
    \item $\mathcal{E}$ is \textit{$k$-entangled} if for each  uncountable family $\mathcal{A}\subseteq [\mathcal{E}]^{k}$ of pairwise disjoint sets and  $t:k\longrightarrow\{<,>\}$ there are distinct $a,b\in \mathcal{A}$ for which that $T(a,b)=t.$
    \item$\mathcal{E}$ is \textit{entangled} if it is $k$-entangled for each $k\in\omega$.
\end{itemize}
\end{definition}
\begin{lemma} Let $(X,<)$ be a total order, $k\in\omega$ and $\mathcal{E}\in [X]^{\omega_1}$ injectively enumerated as $\langle r_\alpha\rangle_{\alpha\in \omega_1}$. Then $\mathcal{E}$ is $k$-entangled if and only if for every uncountable family $\mathcal{C}\subseteq [\omega_1]^{k}$ of pairwise disjoint sets and  each $t:k\longrightarrow\{<,>\}$ there are distinct $c,d\in \mathcal{C}$ for which $$r_{c(i)}\,t(i)\,r_{d(i)}$$ 
for each $i<k.$
\begin{proof}As both implications are proved in a completely similar way, we will only show the one from left to right. For this purpose, suppose that $\mathcal{E}$ is $k$-entangled.   Let $\mathcal{C}\subseteq[\omega_1]^k$ be an uncountable family of pairwise disjoint sets and $t:k\longrightarrow \{<,>\}$.  Given $c\in \mathcal{C}$, take $h_c:k\longrightarrow k$ the unique function satisfying that $i<j$ if and only if $r_{c(h(i))}<r_{c(h_c(j))}$ for every $i,j<k$. By refining $\mathcal{C}$, we can suppose without loss of generality that there is $h$ for which $h=h_c$ for all $c\in \mathcal{C}$. The key observation is that $h$ codes the increasing enumeration (with respect to order in $X$) of $\langle r_{c(i)}\rangle_{i<k}$ whenever $c\in \mathcal{C}$. That is, $r_{c(h(0))}<\dots<r_{c(h(k-1))}$. Now, let  $t'=t\circ h$. By the previous observation and since $\mathcal{E}$ is $k$-entangled, there are distinct $c,d\in \mathcal{C}$ for which $$T(\langle\, r_{c(h(i))}\,\rangle_{i<k},\langle\, r_{d(h(i))}\,\rangle_{i<k})=t'.$$ To finish, take an arbitrary $i<k$ and let $i'=h^{-1}(i)$. Then $r_{c(h(i'))}\,t'(i')\, r_{d(h(i'))}$. But $r_{c(h(i'))}=r_{c(i)}$, $t'(i')=t(i)$ and $r_{d(h(i'))}=r_{d(i)}$. In this way, $r_{c(i)}\,t(i)\, r_{d(i)}$. 
\end{proof}
\end{lemma}

\begin{definition}\label{definitionlex} Let $(X,<)$ be a linear order. 
Given functions $f,g\in X^{\omega}$, we say that $f<_{lex} g$ if $f(n)<g(n)$ where $n=\min(\,k\in\omega\,:\,f(k)\not=g(k)\,)$.
\end{definition}
\begin{rem}\label{remarklex} $(X^\omega,<_{lex})$ is a linear order which can be embedded in $\mathbb{R}$ whenever $X$ is countable.
\end{rem}

\begin{theorem}[Under FCA$^\Delta$]\label{entangledscheme}There is an entangled set.
\begin{proof}Let $\mathcal{F}$ be a fully $\Delta$-capturing scheme of type $\langle m_k,n_{k+1},r_{k+1}\rangle_{k\in\omega}$ 
satisfying that $n_{k+1}\geq 2^{m_k}+1$ for each $k\in\omega.$\\ 
Given $k\in\omega$, let us enumerate $\mathscr{P}(m_{k}\backslash r_{k+1})$, possibly with repetitions, as $\langle C^k_i\rangle_{0<i<n_{k+1}}$ Now, for each $\alpha\in \omega_1$ let $f_\alpha:\omega\longrightarrow \mathbb{Z}$ be given as:
$$f_\alpha(k)=\begin{cases}0 &\textit{if }k=0\textit{ or } \Xi_\alpha(k)=-1\\
\Xi_\alpha(k)&\textit{if }k>0,\,\Xi_\alpha(k)\geq 0 \textit{ and }\lVert\alpha\rVert_{k-1}\in C^{k-1}_{\Xi_\alpha(k)}\\
-\Xi_\alpha(k)&\textit{if }k>0,\,\Xi_\alpha(k)\geq 0\textit{ and }\rVert\alpha\lVert_{k-1}\,\notin C^{k-1}_{\Xi_\alpha(k)}\\
\end{cases}$$

Let $\mathcal{E}=\langle f_\alpha\rangle_{\alpha\in\omega_1}$. We claim that $\mathcal{E}$ is an entangled set in $(\mathbb{Z}^\omega,<_{lex})$. Indeed, let $k\in\omega$,  $t:k\longrightarrow \{>,<\}$  and $\mathcal{C}\subseteq [\omega_1]^k$ be an uncountable 
family of pairwise disjoint sets. As $\mathcal{F}$ is fully $\Delta$-capturing, there is  
$\{c_1,\dots,c_{n_{l}-1}\}\in[\mathcal{C}]^{n_l}$ which is $\Delta$-captured at some level $l>0$. We claim that there is $i<n_l$ so that $T(c_0,c_i)=t$. First of all, note that  $\Delta(c_0(j),c_s(j))$ for each $s<n_l$ and 
$j<k$. This means that $f_{c_0(j)}|_l=f_{c_s(j)}|_l$. Furthermore, since the elements of $\mathcal{C}$ are pairwise disjoint, then $\Xi_l(c_s(j))=\Xi_l(c_s)=s$.  From the two previous facts, it is easy to see that for any $0<s<n_l$ and $j\in\omega$,  the order between $f_{c_0(j)}$ and $f_{c_s(j)}$ is decided by the value of both functions at $l$. Now, let us consider $0<i<n_{l}$ so that $$C^{l-1}_{i}=\{\,\lVert c_0(j)\rVert _{l-1}\,:
\,j<k\textit{ and }t(j)=\,<\,\}.$$ According to the second and third cases of the definitions of the functions, we  have that $0<i=f_{c_i(j)}(l)$ whenever $t(j)=\,<$ and $0>-1=f_{c_i(j)}(l)$ whenever $t(j)=\,>$. This implies $T(c_0,c_{i})=t,$ so we are done.
\end{proof}
\end{theorem}

In \cite{hrusakzindulka}, Hru\v{s}\'ak and Zindulka investigated the (sometimes trivial) ideal $Mon(X)$ of all $\sigma$-monotone subspaces of a given metric space $(X,d)$. Particularly, they showed that every separable metric space of size less than the cardinal invariant $\mathfrak{m}_{\sigma-linked}$ is $\sigma$-monotone. In other words, for all such spaces $X$, the ideal $Mon(X)$ is trivial.  In Question 6.7 of that same paper, they asked: \begin{center}Is there a metric space of cardinality $\omega_1$ that is not $\sigma$-monotone?
\end{center}
We will finish this section by proving that $\text{CA}^\delta$ implies the existence of a metric space of cardinality $\omega_1$ which has no uncountable monotone subspaces. Therefore, it is consistent with arbitrarily large values of the continuum that the previous question has an affirmative answer. 
\begin{rem}\label{remarkdiameter}Let $r, c>0$ and $(Y,d)$ be a $c$-monotone space. Then $(Y,r\cdot d)$ is also $c$-monotone.
\end{rem}

\begin{rem}\label{remarksubspacemonotone}Let $c>0$ and $(Y,d)$ be a $c$-monotone space. If $Z\subseteq Y$, then $Z$ is also $c$-monotone.
\end{rem}

\begin{rem}\label{remarkfinalmonotone}Let $c>0$ and $(Y,d)$ be a $c$-monotone space. Then $Y$ is $c'$-monotone for each $c'>c$.    
\end{rem}

By the previous remarks, it is easy to see that a space $Y$ is not monotone if and only if it is not $\frac{1}{n}$ for each $0<n\in \omega$. The following result is implicit in Lemma 4.4 and Proposition 4.5 of \cite{monotonemetricspaces}.
\begin{lemma}\label{lemmamonotone}Let $n\in\omega$. There is a finite metric space $(Z_n,d)$ so that $Z_n$ is not $\frac{1}{n}$-monotone.
\end{lemma}

Note that if $m$ is a natural number for which there is a space $Z_n$ of size $m$ which is not $\frac{1}{n}$-monotone, then the same is true for each $m'\geq m$.

\begin{theorem}[Under $\text{CA}^\Delta$]\label{nonmonotonescheme}There is a metric $d$ over $\omega_1$ so that $(\omega_1,d)$ has no uncountable monotone subspaces.
\begin{proof}
Let $\tau=\langle m_k,n_{k+1},r_{k+1}\rangle_{k\in\omega} $ be a type so that for any $k>0$, there is a metric $d_k$ over $n_k$ for which the following conditions holds:
\begin{center}For each $0<k'\leq k$, $n_{k'}$ is not a  $\frac{1}{k'}$-monotone subspace of $(n_k,d_k)$.\end{center}  Such type exists due to the Lemma \ref{lemmamonotone} and previous remarks. Furthermore, by Remark \ref{remarkdiameter} we may assume that $diam(n_k)=1$ for each $k$.\\
Given $k\in\omega$, define $s_k=\prod\limits_{i<k}\min(d_{i}(x,y)\,:\,x,y\in n_{i}\text{ and } x\not=y)$ and let $\mathcal{F}$ be a capturing construction scheme of type $\tau$. We define a metric over $\omega_1$ as follows: $$d(\alpha,\beta)=\begin{cases}s_{\Delta(\alpha,\beta)}\cdot d_{\Delta(\alpha,\beta)}(\Xi_\alpha(\Delta(\alpha,\beta)),\Xi_\beta(\Delta(\alpha,\beta)))&\text{ if }\alpha\not=\beta\\
0&\text{ if }\alpha=\beta
\end{cases}
$$
In order to prove that $(\omega_1,d)$ satisfies the conclusion of the theorem, let $X\in[\omega_1]^{\omega}$ and $k\in\omega$. We will prove that  $X$ is not $\frac{1}{k}$-monotone. As $\mathcal{F}$ is $\Delta$-capturing, we can find $D\in[X]^{n_k}$ which is $\Delta$-captured at some level $l>k$. By definition of $d$, note that $d(D(i),D(j))=s_l\cdot d_l(\Xi_{D(i)}(l),\Xi_{D(j)}(l))=s_l\cdot d_l(i,j)$ for any two $i,j<n_k$. In other words, $D$ is isometric to $n_k$ seen as a subspace of $(l,s_l\cdot d_l)$. Therefore, $D$  (an hence $X$) is not $\frac{1}{k}$-monotone.
 \end{proof}   

 \end{theorem}
\section{The basics on the forcing $\mathbb{P}(\mathcal{F})$}\label{sectionforcing}

Let  $\tau=\langle m_k,n_{k+1},r_{k+1}\rangle_{k\in \omega}$ be a (good) type.
 Recall that given $k\in\omega$ and $Y$ a finite set of ordinals of size $m_k$, there is a unique scheme over $Y$ of type $\tau$ which we call $\mathcal{F}(Y)$.  There is also a unique construction scheme of type $\tau$ over $\omega$. Such scheme, denoted by $\mathcal{F}(\omega)$, can be defined as $$\mathcal{F}(
\omega):=\bigcup_{k\in \omega}\mathcal{F}(m_k).$$
In general, it can be proved that for an infinite  $\lambda\leq \omega_1$, there is a construction scheme of domain $\lambda$ if and only if $\lambda$ is limit\footnote{ Such scheme is never unique when $\lambda>\omega$.}. Given a limit $\gamma\leq \lambda$, we define $\mathcal{F}|_\gamma$ as $\{F\in \mathcal{F}\,:\,F\subseteq \gamma\}.$ This family is always a construction scheme with domain $\gamma$. Furthermore, if $\mathcal{F}_\gamma$ is a scheme with domain $\gamma$ so that $\mathcal{F}_\gamma\subseteq \mathcal{F}$, then $\mathcal{F}_\gamma=\mathcal{F}|_\gamma$.\\
According to the previoius paragraph we get that every scheme over $\omega_1$ is fully determined by a sequence of scheme $\langle \mathcal{F}_\gamma\rangle_{\gamma\in \Lim}$ so that:\begin{itemize}
    \item $\mathcal{F}_\gamma$ has domain $\gamma$.
    \item If $\gamma<\lambda$, then $\mathcal{F}_\gamma\subseteq \mathcal{F}_\lambda$.
    \item If $\lambda$ is a limit of limit ordinals, then $\mathcal{F}_\lambda=\bigcup_{\gamma<\lambda}\mathcal{F}_\gamma.$
\end{itemize}
The main difficulty while building such a sequence, is knowing how to extend a scheme $\mathcal{F}$ over a limit $\gamma$ to a new scheme over $\gamma+\omega$. In order to deal with this situation, we use the forcing $\mathbb{P}(\mathcal{F})$ defined below.
This forcing was first introduced by Todor\v{c}ev\'ic in \cite{schemenonseparablestructures} as a tool for building schemes in a recursive manner. Here, we follow the presentation from \cite{schemescruz}. The reader can find the proofs for the results stated in this section in either of those papers.
For the rest of this section, let us fix a construction scheme $\mathcal{F}$ of type $\tau$ over a limit ordinal $\gamma$. 

\begin{definition}[The forcing $\mathbb{P}(\mathcal{F})$]\label{forcingPFdef} We define $\mathbb{P}(\mathcal{F})$ as the forcing consisting of the empty set and  of all $p\in \Fin(\gamma+\omega)$ with the following properties:
\begin{enumerate}[label=(\Roman*), itemsep=0.5em]
\item There is $k_p\in \omega$ such that $|p|=m_{k_p}$.
\item There is $F\in \mathcal{F}_{k_p}$ such that $p\cap \gamma\sqsubseteq F$.
\item $p\cap [\gamma,\gamma+\omega)$ is an initial segment of $[\gamma,\gamma+\omega)$.
\end{enumerate}
Whenever $p\cap\gamma\not=\emptyset$ (even if $p$ is not a condition of $\mathbb{P}(\mathcal{F})$), we let $\alpha_p=\max(p\cap \gamma)$. Additionally, for each $k\in\omega$ we let $\mathbb{P}_k(\mathcal{F})=\{\,p\in\mathbb{P}(\mathcal{F})\,:\,k_p=k\, \}.$  The order on $\mathbb{P}(\mathcal{F})$ is given by $$p\leq q\text{ if and only if }q\in \mathcal{F}(p)\textit{ or }q=\emptyset.$$
\end{definition}
\begin{definition}[The reduction operation]If $p\in \Fin(\gamma+\omega)$ and $\delta\leq \gamma$, we define the \textit{reduction of $p$ to $\delta$} as follows: $$red_\delta(p)=\begin{cases}(p\cap \delta)\cup[\max(p\cap \delta)+1,\max(p\cap \delta)+|p\backslash \delta|+1)&\textit{ if }p\cap\delta\not=\emptyset\\
|p|&\textit{ if }p\cap\delta=\emptyset
\end{cases}
$$
\end{definition}
 Given a filter $\mathcal{G}$ over $\mathbb{P}(\mathcal{F})$, we define  $\mathcal{F}^\mathcal{G}$ as $\bigcup_{p\in \mathcal{G}}\mathcal{F}(p).$ Finally, $\mathcal{F}^{Gen}$ denotes the name for $\mathcal{F}^\mathcal{G}$ where $G$ is a generic filter.

\begin{lemma}\label{equivalencereductioncondition} If $p\in \text{Fin}(\gamma+\omega)$ is such that $p\cap[\gamma,\gamma+\omega)$ is an initial segment of $[\gamma,\gamma+\omega)$, then the following conditions are equivalent:
\begin{enumerate}[label=$(\alph*)$,itemsep=0.5em]
\item $p\in \mathbb{P}(\mathcal{F})$,
\item $red_\gamma(p)\in \mathcal{F}.$
\end{enumerate}
\end{lemma}

\begin{definition}[The cut operation] Let $F\in \Fin(\gamma)$ and $\alpha\in \gamma$. We define the \textit{cut of $F$ at $\alpha$} as follows: $$Cut_\alpha(F)=(F\cap \alpha)\cup [\gamma,\gamma+|F\backslash\alpha|)$$
\end{definition}
The reduction and cut operations are in some sense inverses of each other. That is, if $p\in \mathbb{P}(\mathcal{F})$ and $p\cap \gamma\not=\emptyset$, then 
$Cut_{\alpha_p+1}(red_\gamma(p))=p$. On the other hand, if $F\in \Fin(\gamma)$ and $\alpha\in F$ then $red_\gamma(Cut_{\alpha+1}(F))=F$.

\begin{lemma}\label{cutconditionlemma}Let $F\in \mathcal{F}$ and $\alpha\in \gamma$. Then $Cut_\alpha(F)\in \mathbb{P}(\mathcal{F})$. Furthermore, if $F,G\in \mathcal{F}$ are such that $F\subseteq G$ and $\alpha\in F$, then $Cut_\alpha(G)\leq Cut_\alpha(F)$.

\end{lemma}

\begin{lemma}\label{lemmacut}Let $k\in\omega$, $p\in \mathbb{P}_k(\mathcal{F})$ (that is, $k=k_p$) and $\alpha\in \gamma$ be such that $(\alpha)^-_k=p\cap\gamma$. Then $Cut_\alpha(G)\leq p$ for each $G\in \mathcal{F}$ with $\alpha\in G$ and $\rho^G\geq k.$
\end{lemma}
We now define the property $IH_1$. As it is stated in Proposition \ref{IH1Fomega} below, this property allow us to extend a construction scheme with countable limit domain $\gamma$ to one whose  domain is $\gamma+\omega$.
\begin{definition}[The $IH_1$ property]Let $A\in \Fin(\gamma)$,  $\alpha\in \gamma$ and $F\in \mathcal{F}$. We say that  $IH_1(\alpha,A,F)$ holds if:
\begin{enumerate}[label=$(\arabic*)$,itemsep=0.5em]
    \item $A\subseteq F_0,$
    \item $R(F)=F\cap \alpha$.
\end{enumerate}
Additionally, we say that $\mathcal{F}$ satisfies $IH_1$ if for all $A\in\Fin(\gamma)$ and $\alpha\in \gamma$, there is $F\in \mathcal{F}$ for which $IH_1(\alpha,A,F)$ holds.
\end{definition}

One can prove that $\mathcal{F}(\omega)$ satisfies $IH_1$. Even more, if $\mathcal{F}$ is a scheme with $dom(\mathcal{G})=\lambda$ where $\lambda$ is a limit of limit ordinals, then $\mathcal{F}$ satisfies $IH_1$ if and only if $\mathcal{F}|_\gamma$ satisfies $IH_1$ for each limit  $\gamma<\lambda.$

\begin{proposition}\label{IH1Fomega}Suppose that $\tau$ is a good type. Then $\mathcal{F}(\omega)$ satisfies $IH_1$. Even more, if $\mathcal{F}$ is a construction scheme which satisfies $IH_1$, then there is a countable family  $\mathcal{D}$ of dense sets in $\mathbb{P}(\mathcal{F})$ so that for any filter $\mathcal{G}$ intersecting element of $
\mathcal{D}$,  $\mathcal{F}^\mathcal{G}$ is a construction scheme over $\gamma+\omega$ which  also satisfies $IH_1$.
    
\end{proposition}
\section{$\rho$-capturing and the $\varclubsuit$ principle}\label{sectionclub}
A $C$-sequence over $\omega_1$ is a sequence $\mathcal{C}=\langle C_\delta\rangle_{\delta\in\Lim}$ of subsets of $\omega_1$ of order type $\omega$ with $\sup(C_\delta)=\delta$ for each $\delta\in \Lim$. In the last 30 years, $\mathcal{C}$-sequences have been proved to be a fundamental tool for solving mayor problems in set theory and topology. The reader can consult \cite{Walksonordinals} for more information on $C$-sequences both over $\omega_1$ and  over higher cardinals. The goal of this section is to prove that $\varclubsuit$-principle implies $CA^\rho$. For that, the only thing we need to know is that if $\varclubsuit$-principle holds, then there is a $C$-sequence which is more over a $\varclubsuit$-sequence. \\In what follows, we will fix a  $\varclubsuit$ $C$-sequence $\mathcal{C}=\langle C_\delta\rangle_{\delta\in\text{Lim}}$ and a type $\tau=\langle m_k,n_{k+1},r_{k+1}\rangle_{k\in\omega}$. Note that appart from Theorem \ref{teocapprin}, the analysis presented here works even if $\mathcal{C}$ is not assumed to be $\varclubsuit$ sequence.\\
Let us define $\Lim_0:=\Lim$ and $\Lim_{n+1}:=\{\delta\in \Lim\,:\,C_\delta\subseteq \Lim_n\}$ for each $n\in \omega$. Given $\alpha\in \text{Lim}_0$, put $C^0_\alpha:=C_\alpha$. Finally, for each $n\in\omega$ and $\delta\in \text{Lim}_{n+1}$,  let $$C^{n+1}_\delta:=\bigcup\limits_{\xi\in C_\delta} C^n_\xi.$$
While simple, the next lemma is a fundamental for constructing $\rho$-capturing schemes.
\begin{lemma}\label{lemmatCnstationary}Let $n\in\omega$ and $S\in[\omega_1]^{\omega_1}$.Then  the set $\{\alpha\in \Lim_n\,:\,C^n_\delta\subseteq S\}$ is stationary.
\begin{proof}The proof is carried by induction over $n$. If $n=0$, the conclusion is clear as $\langle C_\delta\rangle_{\delta\in\text{Lim}}$ is a $\varclubsuit$ sequence. Now, suppose that we are given $n\in\omega$ and we already know that $S_n:=\{\delta\in\text{Lim}_n\,:\,C^n_\delta\subseteq S\}$ is stationary. In particular, such set is uncountable, so the set $S_{n+1}:=\{\alpha\in \text{Lim}\,:\,C_\alpha\subseteq S_n\}$ is stationary as well. In order to finish, just note that $S_{n+1}=\{\alpha\in \text{Lim}_{n+1}\,:\,C^{n+1}_\alpha\subseteq S\}.$    
\end{proof}
    
\end{lemma}
The key concepts needed in order to achieve the goal of this section are presented below.
\begin{definition}[Good tuples] Let $\mathcal{F}$ be a construction scheme over an ordinal $\gamma$. We say that a tuple $(n,\delta,k,A)$ is \emph{$\mathcal{F}$-good} if:\begin{itemize}[itemsep=0.5em]
\item$n,k\in\omega$,
\item $\delta\in \text{Lim}_n\cap \gamma$, 
\item $\{ (\xi)_k\,:\,\xi\in C^n_\delta\}$ is a root-tail-tail $\Delta$-system with root $R^n_k(\delta)$.
\item $A\in\text{Fin}(\,[\delta,\gamma)\,)$
\end{itemize}
\end{definition}    

\begin{definition}[The recursive hypothesis] We will say that a construction scheme $\mathcal{F}$ over an ordinal $\gamma$ satisfies $IH_\rho(\mathcal{C},n)$ if for any $(n,\delta,,k,A)$ $\mathcal{F}$-good tuple, there are infinitely many $k\leq l\in \omega$ with $n_l>n$ for which there is $F\in \mathcal{F}_l$ and $D\in [C^n_\delta]^n$ with the following properties:
\begin{enumerate}[itemsep=0.5em]
    \item $A\subseteq F_n\backslash R(F)$.
    \item For each $i<n$, $(D(i))_k\subseteq F_i$ and $(D(i))_k\cap R(F)=R^n_k(\delta)$.
\end{enumerate}
Whenever a pair $(l,F)$ satisfies these properties, we will say that $(l,F)$ accepts $(n,\delta,k,A)$.
\end{definition}
\begin{remark}$\mathcal{F(\omega)}$ satisfies $IH_\rho(\mathcal{C},n)$ for any $n\in\omega$.
    
\end{remark}

\begin{lemma}\label{lemmaballunbounded}Let $\mathcal{F}$ be a construction scheme over a limit ordinal $\gamma$. Also, consider  a limit ordinal $\delta\leq \gamma$, $\alpha\in \delta$ and $k\in\omega$. If $\mathcal{F}$ satisfies $IH_1$, then the set $\{\xi<\delta\,:\,(\xi)^-_k=(\alpha)^-_k\}$ is unbounded in $\delta$.
\begin{proof}If $k=0$ then $(\xi)^-_k=\emptyset$ and the lemma holds trivially. Now, suppose that the lemma holds for some $k\in \omega$. We will prove that it holds for $k+1$. For this purpose, let $\beta\in \delta$. We need to find $\beta\leq\xi<\delta$ so that $(\xi)^{-}_{k+1}=(\alpha)^-_{k+1}$. According to the inductive hypothesis, there is $\beta'\geq \beta$ so that $(\beta')^-_k=(\alpha)^-_k$. In particular, this means that $\rho(\alpha,\beta')>k$.Now, since $\mathcal{F}$ satisfies $IH_1$, there is $F\in \mathcal{F}$ so that $A=\{\alpha,\beta'\}\subseteq F_0$ and $\alpha\cap F=R(F)$.

\end{proof}
\end{lemma}

\begin{lemma}\label{lemaroprin}Let $\mathcal{F}$ be a construction scheme over a countable limit ordinal $\gamma$ and let $n\in\omega$ be such that $n_k>n+1$ for all but finitely many $1\leq k\in\omega$. Assume that $\mathcal{F}$ satisfies $IH_1$, $IH_\rho(\mathcal{C},n)$ and $IH_\rho(\mathcal{C},n+1)$. Then:
$$\mathbb{P}(\mathcal{F})\Vdash\text{\say{$\mathcal{F}^{\text{Gen}}$ satisfies $IH_\rho(\mathcal{C},n+1)$}.}$$
\begin{proof}Let us fix a condition $p\in \mathbb{P}(\mathcal{F})$, a limit ordinal $\delta<\gamma+\omega$, $k\in\omega$ and $A\in\text{Fin}(\, [\delta,\gamma+\omega)\,)$ so that $p\Vdash \text{\say{$(n+1,\delta,k,A)$ is $\mathcal{F}^{\text{Gen}}$-good}.}$ Given $ k'\geq k$, our goal is to find  $l\geq k'$, $q\leq p$ and $F\in \text{Fin}(\gamma+\omega)$ so that $q\Vdash\text{\say{$(l,\mathcal{F})$ accepts $(n,\delta,k,A)$}.}$
Without any loss of generality we may assume that $A\subseteq p$, $k_p\geq k'$ and $p\cap \delta\not=\emptyset$. We will divide the rest of the proof into two cases.\\

\noindent
\underline{Case 1}: $\delta=\gamma$.
\begin{claimproof}[Proof of Case]Let  $\delta'\in C_\gamma$ be such that $\delta'>p\cap\gamma$. By the hypotheses, the set $\{(\xi)_k\,:\,\xi\in C^{n+1}_\gamma\}$ is a root-tail-tail $\Delta$-
system with root $R^{n+1}_k(\gamma)$. Therefore, we may find $\nu\in C^{n+1}_\gamma$ such that $\min(\,(\nu)_k\backslash R^{n+1}_k(\gamma)\,)>\delta'$. Now,  let us consider $\alpha'=\min (\text{red}_\gamma(p)\backslash p)$. Note that $(\alpha)^-_{k'}=\alpha\cap p=\gamma\cap p$. Therefore, we can make use of  Lemma \ref{lemmaballunbounded} to find $\nu<\alpha<\gamma$ such that $(\alpha)^-_{k'}=\gamma\cap p$. We now define $$A':=\big( (\nu)_k\backslash R^{n+1}_k(\gamma)\big)\cup \{\alpha\}.$$
Since $\gamma\in\text{Lim}_{n+1}$, then $\delta'\in \text{Lim}_n$. Furthermore, $C^n_{\delta'}\subseteq C^{n+1}_{\gamma+1}$. Thus, the family  $\{\,(\xi)_k\,:\,\xi\in C^n_{\delta'}\}$ is also a root-tail-tail $\Delta$-system with root $R^n_k(\delta')=R^{n+1}_k(\gamma)$. In this way, we conclude that $(n,\delta',k, A')$ is a $\mathcal{F}$-good tuple. Recall that $\mathcal{F}$ is assumed to satisfy $IH_\rho(\mathcal{C},n)$.  
Thus, we may find $l>k'$ with $n_l>n+1$ and $F\in \mathcal{F}_l$ such that $(l,F)$ accepts $(n,\delta',k,A')$. In particular, this means that $\alpha\in F_n\backslash{R(F)}$. By the homogeneity, there is $\beta\in F_{n+1}\backslash {R(F)}$ such that $(\beta)^-_{k'}=(\alpha)^-_{k'}$. We now can define $q:=\text{Cut}_\beta(F)$. Note that $q\leq \text{Cut}_\alpha(F_{n+1})\leq p$. It is easy to see that $$q\Vdash\text{\say{$(l,q)$ accepts $(n+1,\gamma,k,A)$}.}$$
\end{claimproof}
\noindent
\underline{Case 2}: $\delta<\gamma$.
\begin{claimproof}In this case,  we first use Lemma \ref{lemmaballunbounded} to get  $\max(\,(A\cap \gamma)\cup\{\delta\}\,)<\alpha<\gamma$ with $(\alpha)^-_{k'}=\gamma\cap p$. We now define $A'=(A\cap \gamma)\cup \{\alpha\}.$
Since $\delta \in \text{Lim}_{n+1}\cap \gamma$ and $\mathcal{F}$ satisfies $IH_\rho(\mathcal{C},n+1)$, we can find $l>k'$ and $F\in \mathcal{F}$ so that the pair $(l,F)$ accepts $(n+1,\delta,k,A')$. Let us define $q:=\text{Cut}_\alpha(F)$. Then $q\leq Cut_\alpha(F_{n+1})\leq p$. As in the previous case, we have that $$q\Vdash\text{\say{$(l,q)$ accepts $(n+1,\gamma,k,A)$}.}$$ \end{claimproof}
\end{proof}
\end{lemma}
 
The  property $IH_\rho(\mathcal{C},n+1)$ can be coded by countably many dense sets in the forcing $\mathbb{P}(\mathcal{F})$. Thus, Lemma \ref{lemaroprin} can be rewritten as:
\begin{corollary}Let $\mathcal{F}_\gamma$ be a construction scheme over some countable ordinal $\gamma$  and let $n\in\omega$ be such that $n_k>n+1$ for each $1\leq k\in \omega$. Assume that $\mathcal{F}_\gamma$ satisfies $IH_1$, $IH_\rho(\mathcal{C},n)$ and $IH_\rho(\mathcal{C},n+1)$. Then there is a construction scheme $\mathcal{F}_{\gamma+\omega}$ which also satisfies $IH_1$, $IH_\rho(\mathcal{C},n)$ and $IH_\rho(\mathcal{C},n+1)$ and such that $\mathcal{F}_\gamma\subseteq \mathcal{F}_{\gamma+\omega}$.  
\end{corollary}

The lemma above serves as the successor step for a recursive construction. For limit steps, we will make use of the following easy lemma whose proof we omit.
\begin{lemma}Fix $n\in\omega$. Let  $\gamma\leq \omega_1$ be a limit of limit ordinals  and let $\mathcal{F}_\gamma$  be a construction scheme over $\gamma$ such that $\mathcal{F}\gamma|_{\delta}$ satisfies $IH(\mathcal{C},n)$. Then  $\mathcal{F}_\gamma$ also satisfies $IH(\mathcal{C},n)$.
\end{lemma} 
As a direct consequence of the two results above, we have the following.
\begin{theorem}There is a construction scheme over $\omega_1$ which satisfies $IH_1$ and $IH(\mathcal{C},n)$ for each $n\in\omega$ with the property that $n_k>n$ for all but finitely many $k$'s.
\end{theorem}
\begin{theorem}\label{teocapprin}If $\mathcal{F}$ is a construction scheme over $\omega_1$ which satisfies $IH_{\rho}(\mathcal{C},n)$, then $\mathcal{F}$ is $n+1$-$\rho$-capturing.
\begin{proof}Let $\mathcal{S}$ be an uncountable subset of $\text{Fin}(\omega_1)$. By the $\Delta$-system lemma, there is a root-tail-tail $\Delta$-system $\mathcal{S}_1\in[\mathcal{S}]^{\omega_1}$ with some root $R$. By refining $\mathcal{S}_1$ we may assume that there is $k\in\omega$ so that $\rho^D=k$ for any $D\in S_1$. Now, let $\alpha_D:=\max(D)$ for any $D\in \mathcal{S}_1$. Note that $D\subseteq (\alpha_D)_k$. Again, according to the $\Delta$-system Lemma, there is $\mathcal{S}_2\in[\mathcal{S}_1]^{\omega_1}$ such that the family $\{(\alpha_D)_k\,:\,D\in \mathcal{S}_2\}$ forms a root tail-tail $\Delta$-system with some root $R'$ with the extra property that $R'\cap D=R$ for each $D\in \mathcal{S}_2$ (In other words, $D\backslash R\subseteq (\alpha_D)_k\backslash R')$.

By virtue of Lemma \ref{lemmatCnstationary}, we can find $\delta\in \text{Lim}_n$ for which $C^n_\delta\subseteq \{\alpha_D\,:\,D\in\mathcal{S}_2\}$. In order to finish, let $E\in \mathcal{D}$ for which $$A:=(\alpha_E)_k\backslash R'\subseteq [\delta,\omega_1).$$ Note that the  tuple $(n,\delta,k, A) $ is $\mathcal{F}$-good. Since $\mathcal{F}$ satisfies $IH_\rho(\mathcal{C},n)$, there is $l\geq k$ with $n_l>n$, $F\in \mathcal{F}_l$ and $D_0,\dots,D_{n-1}\in [\mathcal{S}_2]^n$  with $\{\alpha_{D_i}\,:\,i<n\}\subseteq C^n_\delta$ such that the following properties hold:
\begin{enumerate}
    \item $A\subseteq F_n\backslash R(F)$.
    \item For each $i<n$, $(\alpha_{D_i})_k\subseteq F_i$ and $(\alpha_{D_i})_k\cap R(F)=R'.$
\end{enumerate}
We claim that $\{D_0,\dots,D_{n-1},E\}$ is  $\rho$-captured at level $l$. For this purpose, we need to show that $D_i\subseteq F_i$ and $D_i\cap R(F)=R$ for each $i<n$, and that $E\subseteq F_n$ and $E\cap R(F)=R$. Indeed, for the first part, note that $D_i\subseteq (\alpha_{D_i})_k\subseteq F_i$ and $$D_i\cap R(F)=D_i\cap (\alpha_{D_i})_k\cap R(F)=D_i\cap R'=R.$$
 Finally, note that $(\alpha_E)=R'\cup A\subseteq F_i$ and $(\alpha_E)_k\cap R(F)=R'.$ Therefore, $E\subseteq F_n$ and $E\cap R(F)=R$ due to the same argument as in the previous paragraph. 
\end{proof}
\end{theorem}
\begin{corollary}[Under $\varclubsuit$-principle] $CA^\rho$ holds.
\end{corollary}
Using slight modifications of the definitions given above, it is easy to se that $CA^\rho(part)$ assuming $\varclubsuit$.

\section{$\Delta$-capturing and $CH$}\label{sectionch}
For the rest of this section we will assume that $CH$ holds. Under this hypothesis we can enumerate $\mathscr{P}(\omega_1)$ as $\mathcal{C}=\langle C_\gamma\rangle_{\gamma\in \text{Lim}}$ in such way that $
\beta_\gamma:=\sup(C_\gamma)<\gamma$ for any $\gamma$. We fix such an enumeration. Finally, we fix a type $\tau=\langle m_k,n_{k+1},r_{k+1}\rangle_{k\in\omega}$.

\begin{definition}Let $\mathcal{F}$ be a construction scheme over a limit ordinal $\gamma$, $\delta\in \text{Lim}\cap \gamma$, $\alpha\in [\delta,\gamma)$, $l\in\omega$ and $D\in \text{Fin}(C_\delta)$.  We say that $D$ is \emph{$(\delta,l,\alpha)$-adequate} if the following properties hold:
\begin{enumerate}
\item $D$ is $\Delta$-captured at level $l$.

    \item For any $i<j$, $\Delta(\alpha,D(i))\geq l$.
\end{enumerate}
We define $j^\alpha_\delta(l):=\max(\,j\leq n_l\,:\,\exists D\in [C_\delta]^j\,(D\text{ is }(\delta,l,\alpha)\text{-adequate})\,).$
\end{definition}

\begin{definition}We will say that a construction scheme over an ordinal $\gamma$ satisfies $IH_\Delta(\mathcal{C})$ if for any $\delta \in \text{Lim}\cap \gamma$ and each $\alpha\in [\delta,\gamma)$  there are infinitely many $l\in \omega$ for which $r_l=|(\alpha)_{l-1}\cap \delta|$ and either $j^\alpha_\delta(l)=n_l$ or there is $\Xi_\alpha(l)=j^\alpha_\delta(l)$.
\end{definition}
\begin{lemma}Let $\mathcal{F}$ be a construction scheme over a countable limit ordinal $\gamma$ which satisfies $IH_1$. Then $\mathbb{P}(\mathcal{F})\Vdash\text{\say{$\mathcal{F}^{\text{Gen}}$ satisfies $IH_\Delta(\mathcal{C})$}}.$
\begin{proof}
Let $p\in\mathbb{P}(\mathcal{F})$, $\delta\in\text{Lim}\cap (\gamma+\omega)$, $\alpha\in [\delta,\gamma+\omega)$ and $k\in\omega$. Assume that there is no $q'\leq p$ and $l,j\in\omega$ for which $q \forces{$j=j^\alpha_\delta(l)=n_l$ }$. Our goal is to find $q\leq p$ and $l,j\in\omega$ such that $ q\forces{$r_l=|(\alpha)_{l-1}\cap \delta|\text{ and }j=j^\alpha_\delta(l)=\Xi_\alpha(l)$}$ . If $\alpha\in \gamma$, we are done. This is because $\mathcal{F}$ is forced to be contained in $\mathcal{F}^{\text{Gen}}$ and $\mathcal{F}$ already satisfies $IH_\Delta(\mathcal{C})$. Thus, we may assume that $\alpha\in [\gamma,\gamma+\omega)$. Without any loss of generality we can also suppose that $\alpha\in p$ and $k_p\geq k$.  We divide rest of the proof into two cases.\\

\noindent
\underline{Case 1}: $\delta=\gamma$.
\begin{claimproof}[Proof of case] Since $\mathcal{F}$ satisfies $IH_1$, we can extend $p$ to a condition $q'$ such that $r_{k_q'}+1=|q'\cap \gamma|$. Here, we let $l=k_{q'}+1$ and $j$ be the maximum of all $j'\leq n_l$ for which there is $D\in [C_\gamma]^{j'}$ which is $\Delta$-captured at level $l$ and such that $\lVert D(i)\rVert_{l-1}=|\alpha\cap q'|$ for each $i<j$. As $\alpha\in q'$ and $k_{q'}=l-1$, then $q'\forces{$\alpha\cap q'=(\alpha)^-_{l-1}$}.$  Therefore, $$q'\forces{$\lVert\alpha\rVert_{l-1}=|\alpha\cap q'|\text{ and }r_l=|(\alpha)_{l-1}\cap \gamma|  $ }.$$ From this fact and by maximality of $j$, it also follows that $q\forces{$j=j^\alpha_\gamma(l)$}$. In order to finish, let $F\in \mathcal{F}_l$ be such that $q\cap \gamma\subseteq F.$ Note that since $r_l=|q'\cap \gamma|$, then $R(F)=q'\cap \gamma$. Now, according to our initial hypotheses, $q'\forces{$j<n_l$}$. In this way, $F_j$ is well-defined. Let $q:=\text{Cut}_{\alpha_j}(F)$ where $\alpha_j=\min(F_j\backslash R(F))$. Then $q\forces{$\Xi_\alpha(l)=j$}$, and since $(\alpha_j)_{l-1}=R(F)=q'\cap \gamma$, then $q\leq q'$. 
\end{claimproof} 

\noindent
\underline{Case 2}: $\delta<\gamma$.
\begin{claimproof}[Proof of case]
    
Without any loss of generality we may assume that $p\cap [\delta,\gamma)\not=\emptyset$. Let $F=\text{red}_\gamma(p)$ and $\alpha'=F(|\alpha\cap p|)$ and $\beta=F(|p\cap \gamma|)=\max(p\cap \gamma)+1$. Note that $\alpha'\geq \beta\geq \delta$. Since $\mathcal{F}$ satisfies $IH_\Delta(\mathcal{C})$ and $\delta\leq \alpha'<\gamma$, then there is $l>k_p$ with $r_l=|(\alpha)_{l-1}\cap \delta|$ and such that either $j^\alpha_\delta(l)=n_l$ or $\Xi_\alpha(l)=j^\alpha_\delta(l)$. Leg $G\in \mathcal{F}_l$ be such that $\alpha\in G$ and let $j=\Xi_\alpha(l)$. Then $G \cap \delta=R(G)$. In this way, $j\geq 0$ and $\beta,\alpha\in F_j\backslash R(F)$. We define $q:=Cut_\beta(G)$. Since $(\beta)_{k_p}=p\cap \gamma$, then $q\leq p$. Furthemore, $q\forces{ $\Delta(\alpha, \alpha')>l$ }$ and . In this way, $$q\forces{$j^\alpha_\delta(l)=j^{\alpha'}_\delta(l)=j $}.$$ In particular, this means that $j<n_l$, so $q\forces{$\Xi_\alpha(l)=\Xi_{\alpha'}(l)=j$  }$. Finally, since $\beta\geq \delta$, then $q\forces{$(\alpha)_{l-1}\cap\delta= (\alpha')_{l-1}\cap \delta$}$. This means that $q\forces{$r_l=|(\alpha)_{l-1}\cap\delta|$}$. Thus, the proof is over.
\end{claimproof}

\end{proof}
    
\end{lemma}
\begin{theorem}Let $\mathcal{F}$ be a construction scheme over $\omega_1$ which satisfies $IH_\Delta(\mathcal{C})$. Then $\mathcal{F}$ is fully $\Delta$-capturing.
\begin{proof}Let $\mathcal{S}$ be an uncountable subset of $\text{Fin}(\omega)$. By the $\Delta$-system Lemma, there is a root-tail-tail $\Delta$-system $\mathcal{S}_1\in [\mathcal{S}]^{\omega_1}$ with some root $R$. By refining $\mathcal{S}_1$ we may assume that there are $k ,a\in\omega$  so that for any $D,E\in \mathcal{S}$:
\begin{enumerate}
\item$\rho^D=k$.
\item$a=\lVert \alpha_D\rVert_k$ where $\alpha_D=\max(D)$.
If $\phi:(\alpha_D)_k\longrightarrow (\alpha_E)_k$ is the increasing bijection, then $\phi[D]=E$.
\end{enumerate}
Let $M$ be a countable elementary submodel of $H(\omega_2)$ be such that $\mathcal{S}_1$ and $\mathcal{F}$ are in $M$. Now, consider $\gamma>\omega_1\cap M$ such that $C_\gamma=\{\alpha_E\,:\,E\in \mathcal{S}_1\cap M\}$. Since $\mathcal{S}_1$ is a root-tail-tail $\Delta$-system, there is $D\in \mathcal{S}_1$ for which $D\backslash R\subseteq [\gamma,\omega_1)$. Note that $D\cap M= R$.
As $\mathcal{F}$ is assumed to satisfy $IH_\Delta(\mathcal{C})$, there is $k\leq l\in \omega$ with $r_l=|(\alpha)_{l-1}\cap \gamma|$ such that either $j=n_l$ or $\Xi_\alpha(l)=$ where $j:=j^{\alpha_D}_\gamma(\alpha)$.\\

\noindent
\underline{Claim}: $j=n_l$.
\begin{claimproof}[Proof of claim]Suppose towards a contradiction that this is not the case. Then $0\leq j<n_l$ and $\Xi_\alpha(l)=j$. By elementarity, there is $D_j\in \mathcal{S}_1\cap M$ such that $\lVert\alpha_D\rVert_l=\lVert \alpha_{D_j}\rVert_l$. Note that $\Delta(\alpha_{D_j},\alpha_D)>l$ Now, according to the definition of $j$ we can find $D_0,\dots,D_{j-1}\in\mathcal{S}_1\cap M$ such that $\{\alpha_{D_0},\dots,\alpha_{D_{j-1}}\}$ is $\Delta$-captured and $\Delta(\alpha_{D_i},\alpha_{D})\geq l$ for each $i<j$. Since $\Xi_{\alpha_D}(l)=j$, then $\{\alpha_{D_0},\dots,\alpha_{D_{j-1}},\alpha_{D}\}$ is also  $\Delta$-captured at the same level. This contradicts the maximality of $j$, so we are done.
\end{claimproof}

Since $j=n_l$, there are $D_0,\dots, D_{n_l-1}$ for which  $\alpha_{D_0},\dots, \alpha_{D_{n_l-1}}$  are fully $\Delta$-captured. By (1) and (2), it easily follows that  $\mathcal{D}=\{D_0,\dots,D_{n_l-1}\}$ is fully $\Delta$-captured as well at the same level.
\end{proof}
\end{theorem}
It is direct that $\mathcal{F}(
\omega)$ satisfies both $IH_1$ and $IH_\Delta(\mathcal{C})$. Furthermore, the property $IH_\Delta(\mathcal{C})$ is preserved under taking unions. In other words,
\begin{lemma}Fix $n\in\omega$. Let  $\gamma\leq \omega_1$ be a limit of limit ordinals  and let $\mathcal{F}_\gamma$  be a construction scheme over $\gamma$ such that $\mathcal{F}\gamma|_{\delta}$ satisfies $IH_\Delta(\mathcal{C})$. Then  $\mathcal{F}_\gamma$ also satisfies $IH_\Delta(\mathcal{C})$.
\end{lemma} 
Thus, arguing in a similar way as we did in Section \ref{sectionclub}, we conclude :

\begin{corollary}[Under $CH$]\label{fcadeltach} $FCA^\Delta$ holds.
\end{corollary}

Using slight modifications of the definitions given above, it is easy to prove that $FCA^\Delta(part)$ follows from CH.
\section{Open questions}\label{sectionfinal}

One of the most famous open problems is whether there is a normal not countably paracompact space of size $\omega_1$, i.e., a Dowker space. Examples of such spaces exist under various set-theoretic assumptions such as the existence of Lusin set  \cite{pauldowker} (first proved in \cite{PartitionProblems}), and more generally the principle $\varclubsuit _{\rm AD}(\mathcal{ S},\mu, \theta)$ in \cite{smalldowkershalev}, some of whose variants follow from both $\varclubsuit$ or the existence of a Suslin tree. It would be interesting if such spaces can be constructed using $\rho$-capturing schemes.
\begin{problem}Does the existence of a Dowker space of size $\omega_1$ follow from $CA^\rho$?
\end{problem}

\begin{problem}Is the existence of a $3$-$\rho$-capturing scheme  consistent with the non-existence of Suslin trees?
\end{problem}
The proof of theorem \ref{luzinjonestheorem} relies only in the fact that both $n$-$\rho$-captured and $n$-$\Delta$-captured families satisfy the condition (I) in Definition \ref{capturedsystemdef}. These types of schemes are what we would call $\emptyset$-capturing schemes. We wonder whether the capturing axiom associated to these type of schemes, namely $CA^\emptyset$ follows from any guessing principle which is a consequence of both $\varclubsuit$ and CH. Concretely:

\begin{problem}Does the axiom $CA^\emptyset$ follow from $\stick$-principle?
\end{problem}

\begin{problem} Does the existence of an $n$-inseparable $AD$ family imply the existence of an $n$-$\omega_1$-gap? What about the other direction?
\end{problem}

Of course, the existence of an $n$-inseparable $AD$ family which codes $n\times\omega_1$ is enough to imply the existence of an $n$-$\omega_1$-gap. Such a family should at least be Jones. 

\begin{problem}Does the existence of an $n$-inseparable $AD$ family imply the existence of one which is moreover Jones?  
\end{problem}

There are two known constructions of an inseparable Jones family in $ZFC$. The one in this paper and the one in \cite{guzman2019mathbb}. Both of them are not only inseparable but according to Definition \ref{Luzinrepdef}, but they both code $\omega_1$.  We then ask:

\begin{problem}Is it true that every inseparable Jones AD family codes $\omega_1$?

\end{problem}

As for the gap cohomology group, specifically Theorem \ref{biggapcohomologytheorem} and the observation regarding strong-$Q$-sequences prior to it, we have the following question.
\begin{problem}Is it consistent that every $|\mathcal{G}(\mathcal{T})|=2^{\omega_1}$ for any $*$-lower semilattice $\mathcal{T}$?
\end{problem}
\bibliographystyle{plain}
\bibliography{bibliografia}
{\color{white}hlj}\\\\
\noindent
Jorge Antonio Cruz Chapital\\
Department of Mathematics, University of Toronto, Canada\\
cruz.chapital at utoronto.ca\\
\end{document}